\documentclass[11pt]{article}

\usepackage{epic,eepic,epsfig,amssymb,amsmath,amsthm,graphics}

\usepackage{float}
\usepackage{authblk}
\usepackage{multimedia}
\usepackage{xcolor}
\usepackage{tikz}
\usepackage{algorithm2e}
\usepackage{pgfplots}
            
\usepackage[margin=1in]{geometry}

\AtEndDocument{\bigskip{\footnotesize%
  \textsc{Laboratoire J.A. Dieudonné, Université Côte d'Azur, Parc Valrose, 06108 Nice Cedex 2, France} \par  
  \textit{E-mail address}: \texttt{ludovic.rifford@math.cnrs.fr} 
}}
  
\def\R{\mathbb R}

\def\S{\mathbb S}

\numberwithin{equation}{section}

\newtheorem{theorem}{Theorem}

\newtheorem{lemma}[theorem]{Lemma}

\newtheorem{remark}[theorem]{Remark}

\numberwithin{theorem}{section}

\makeindex

\begin{document}

\title{A quantitative version of Tao's result on the Toeplitz Square Peg Problem}


\author{Ludovic Rifford}



\maketitle

\begin{abstract}
Building on a result by Tao, we show that a certain type of simple closed curve in the plane given by the union of the graphs of two $1$-Lipschitz functions inscribes a square whose sidelength is bounded from below by a universal constant times the maximum of the difference of the two functions.
\end{abstract}

\section{Introduction}\label{SECintroduction}

A subset $\Gamma$ of the plane $\R^2$ is said to {\it inscribe a square} if it contains the four vertices of a square with positive sidelength. The Square Peg Problem  raised by Toeplitz \cite{toeplitz} in 1911  can be stated as follows (we recall that a (continuous) curve $\gamma:[0,1] \rightarrow \R^2$ is called closed if $\gamma(0)=\gamma(1)$ and simple if the function $t\in [0,1) \mapsto \gamma(t)$ is injective):\\

\noindent {\bf Square Peg Problem.}
{\em Let $\gamma:[0,1] \rightarrow \R^2$ be a simple closed  continuous curve. Does $\gamma ([0,1])$ necessarily inscribe a square?}

\begin{figure}[H]
\begin{center}
\includegraphics[width=4cm, height=4.1cm]{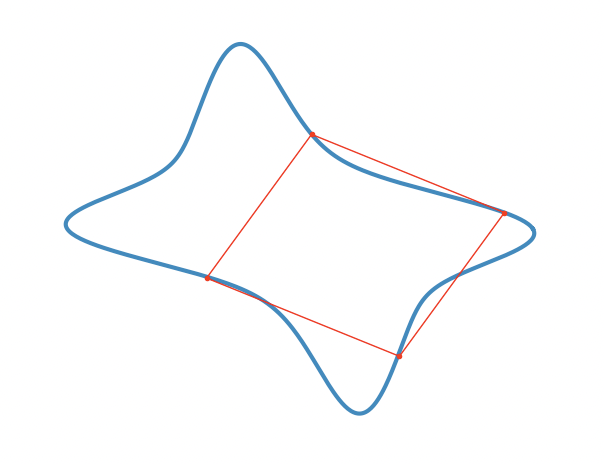}
\caption{A square in red inscribed in the blue curve \label{fig1}}
\end{center}
\end{figure}

The answer to the Square Peg Problem is known to be "Yes" for curves with enough regularity ({\it e.g.} convex, piecewise analytic or locally monotone curves) but remains open in its full generality (in the case of merely continuous simple closed curves). For further details, we refer for example the interested reader to the survey by Matschke \cite{matschke14}. The absence of positive result in the continuous case is due, in particular, to the lack of positive lower bound for the sidelengths of squares inscribed in smooth simple closed curves (or any curve in a set which is dense, in some sense, in the set of continuous simple closed curves). As a matter of fact, if we could show for instance that smooth simple closed curves always inscribe a square whose sidelength is bounded from below by some quantity depending continuously on the curve (such as for example the area enclosed by the curve or its diameter) then it would allow us to prove, by a simple  argument of approximation,  that any (continuous) simple closed curve inscribes a square. The aim of the present paper is precisely to show that we can bound from below the sidelength of squares inscribed in the type of sets investigated by Tao in \cite{tao17}.\\

\begin{figure}[H]
\centering
\begin{tikzpicture}[scale=1]
\draw[color=gray] (0.93,0.768) -- (0.866,-1);
\draw[color=gray] (-0.882,-0.935) -- (0.866,-1);
\draw[color=gray] (-0.882,-0.935) -- (-0.817,0.813);
\draw[color=gray] (0.931,0.748) -- (-0.817,0.813);
\draw[color=gray] (0.8714,-0.8543) -- (0.7257,-0.8488);
\draw[color=gray]  (0.7204,-0.9946) -- (0.7257,-0.8488);
\draw[color=gray]  (0.9985,-0.026) -- (0.7985,-0.226);
\draw[color=gray]  (-0.7495,0.039) -- (-0.9495,-0.161);
\draw[color=gray]  (0.157,0.8805) -- (-0.043,0.6805);
\draw[color=gray]  (0.1,-0.8675) -- (-0.1,-1.0675);
\fill[color=black] (0.866,-1) circle (0.25mm);
\fill[color=black] (-0.882,-0.935) circle (0.25mm);
\fill[color=black] (0.931,0.748) circle (0.25mm);
\fill[color=red] (-0.817,0.813) circle (0.25mm);
\draw (0.866,-1) node[right]{$O$};
\draw (0.931,0.748) node[right]{$P$};
\draw (-0.882,-0.935) node[left]{$R$};
\draw (-0.817,0.813) node[left]{{\color{red} $Q$}};
\end{tikzpicture}
\caption{The point $R$ is the image of $P$ by the rotation of angle $\pi/2$ about $O$ and $Q$ is the opposite corner to $O$ in the square $(OPQR)$\label{fig2}}
\end{figure}
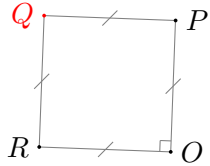

Another way to state the Square Peg Problem is to see it as a problem of intersection of a set with its set of opposite square corners. We say that a triple $(O,P,R)$ in $(\R^2)^3$ is a {\it square corner} if the three points $O,P,R$ are distinct and if
\[
R=\mbox{Rot}_O^{\pi/2}(P),
\] 
where $\mbox{Rot}_O^{\pi/2}:\R^2 \rightarrow \R^2$ denotes the rotation  of angle $\pi/2$ about the point $O$. Then, denoting by $\mathcal{SC}\subset (\R^2)^3$ the set of square corners, we call {\it opposite corner} of a triple $(O,P,R)\in \mathcal{SC}$  the unique point $Q=\mathcal{Q} (O,P,R)$ which makes the quadrilateral $(OPQR)$ a square (see Figure \ref{fig2}), that is 
\[
\mathcal{Q} (O,P,R):= P + \overrightarrow{OR}.
\]

\begin{figure}[H]
\centering
\begin{tikzpicture}[scale=1]
\draw[color=blue] (0,0) ellipse (1cm and 2 cm);
\draw[color=black!40!green] (-0.133,-1.866) ellipse (2cm and 1 cm);
\draw[color=gray] (0.93,0.768) -- (0.866,-1);
\draw[color=gray] (-0.882,-0.935) -- (0.866,-1);
\draw[color=gray] (-0.882,-0.935) -- (-0.817,0.813);
\draw[color=gray] (0.931,0.748) -- (-0.817,0.813);
\draw[color=gray] (0.8714,-0.8543) -- (0.7257,-0.8488);
\draw[color=gray]  (0.7204,-0.9946) -- (0.7257,-0.8488);
\draw[color=gray]  (0.9985,-0.026) -- (0.7985,-0.226);
\draw[color=gray]  (-0.7495,0.039) -- (-0.9495,-0.161);
\draw[color=gray]  (0.157,0.8805) -- (-0.043,0.6805);
\draw[color=gray]  (0.1,-0.8675) -- (-0.1,-1.0675);
\fill[color=black] (0.866,-1) circle (0.25mm);
\fill[color=black] (-0.882,-0.935) circle (0.25mm);
\fill[color=black] (0.931,0.748) circle (0.25mm);
\fill[color=red] (-0.817,0.813) circle (0.25mm);
\draw (0.866,-0.9) node[right]{$\gamma(t)$};
\draw (0.931,0.748) node[right]{$\gamma(u)$};
\draw (-0.882,-0.935) node[left]{$\gamma(v)$};
\draw (-0.817,0.813) node[left]{\color{red} $\mathcal{Q}(\gamma(t), \gamma(u), \gamma(v)) $};
\draw (1.95,-1.866) node[right]{\color{black!40!green} $\mbox{Rot}_{\gamma(t)}^{\pi/2} (\Gamma)$};
\end{tikzpicture}
\caption{The set $\mathcal{SOSC}(\Gamma)$  is the union over $t\in [0,1)$ of opposite corners of the form $\mathcal{Q}(\gamma(t),\gamma(u),\gamma(v))$ where $\gamma(u), \gamma(v)$ are such that $(\gamma(t),\gamma(u),\gamma(v))$ is a square corner \label{fig3}}
\end{figure}
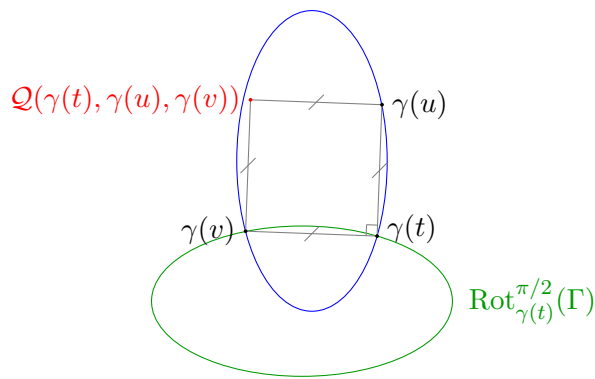

Now, given a subset $\Gamma$ of the plane, we define its {\it set of opposite square corners}, denoted by $\mathcal{SOSC}(\Gamma)$, as the set of opposite corners of all square corners in $\Gamma$, that is 
\[
\mathcal{SOSC}(\Gamma) := \Bigl\{ \mathcal{Q} (O,P,R) \, \vert \, (O,P,Q) \in \mathcal{SC}\cap \Gamma^3 \Bigr\}.
\]
If $\Gamma=\gamma([0,1])$ with $\gamma:[0,1] \rightarrow \R^2$  a simple closed curve, then this set can also be written as  (see Figure \ref{fig3})
\[
\mathcal{SOSC}(\Gamma) := \bigcup_{t,u,v\in [0,1)} \Bigl\{ \mathcal{Q}(\gamma(t), \gamma(u), \gamma(v)) \in \R^2 \, \vert \, u \neq t  \mbox{ and } \gamma(v) = \mbox{Rot}_{\gamma(t)}^{\pi/2} (\gamma(u)) \Bigr\}.
\] 

By construction of the set of opposite square corners, the Square Peg Problem is equivalent to asking whether the set $\Gamma \cap \mathcal{SOSC}(\Gamma)$ is empty or not:  A set $\Gamma = \gamma ([0,1])$,  with $\gamma:[0,1] \rightarrow \R^2$  a simple closed curve, does inscribe a square if and only if the set $\Gamma \cap \mathcal{SOSC}(\Gamma)$ is not empty.\\

\begin{figure}[H]
\begin{center}
\includegraphics[width=4cm,  height=4cm]{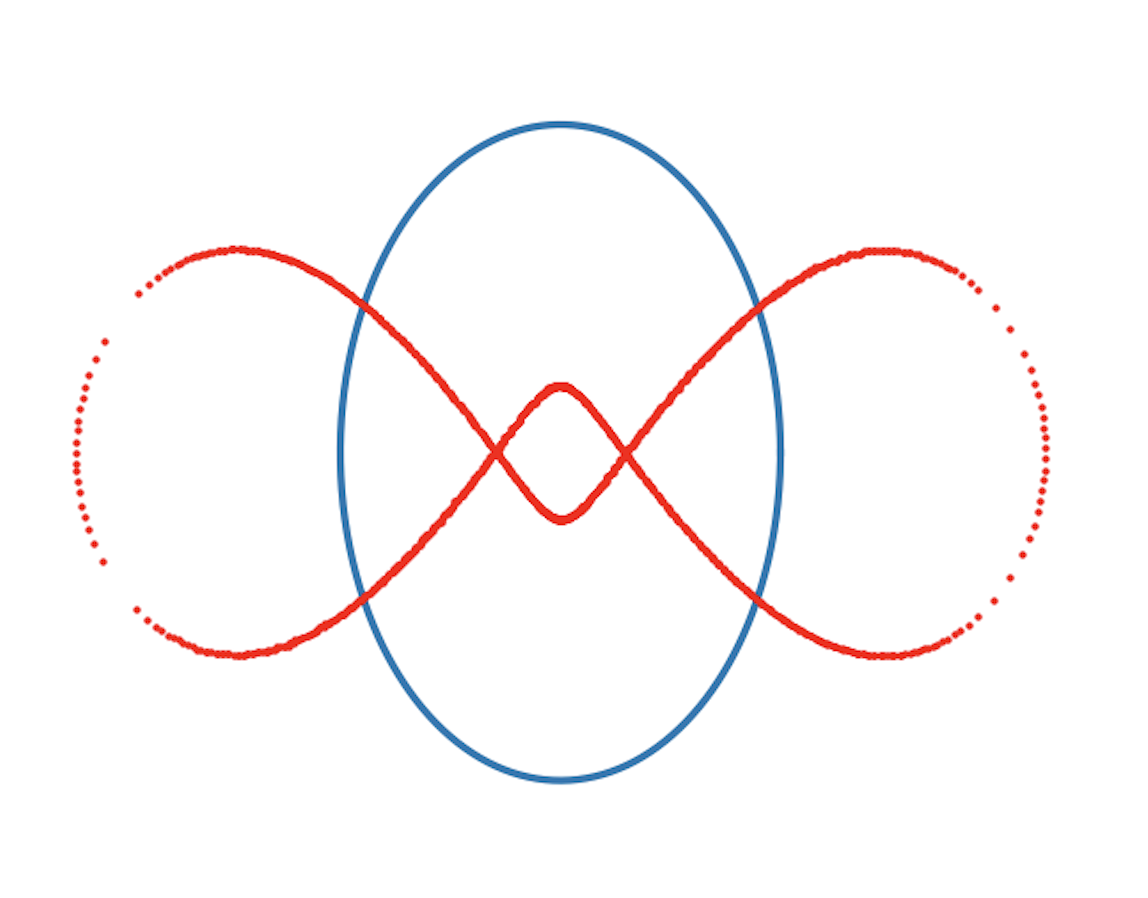}
\caption{In red a discretization of $\mathcal{SOSC}(\Gamma)$ with $\Gamma$ the ellipse of equation $4x^2+y^2=4$ \label{fig4}}
\end{center}
\end{figure}

\begin{figure}[H]
\begin{center}
\includegraphics[width=4cm, height=4cm]{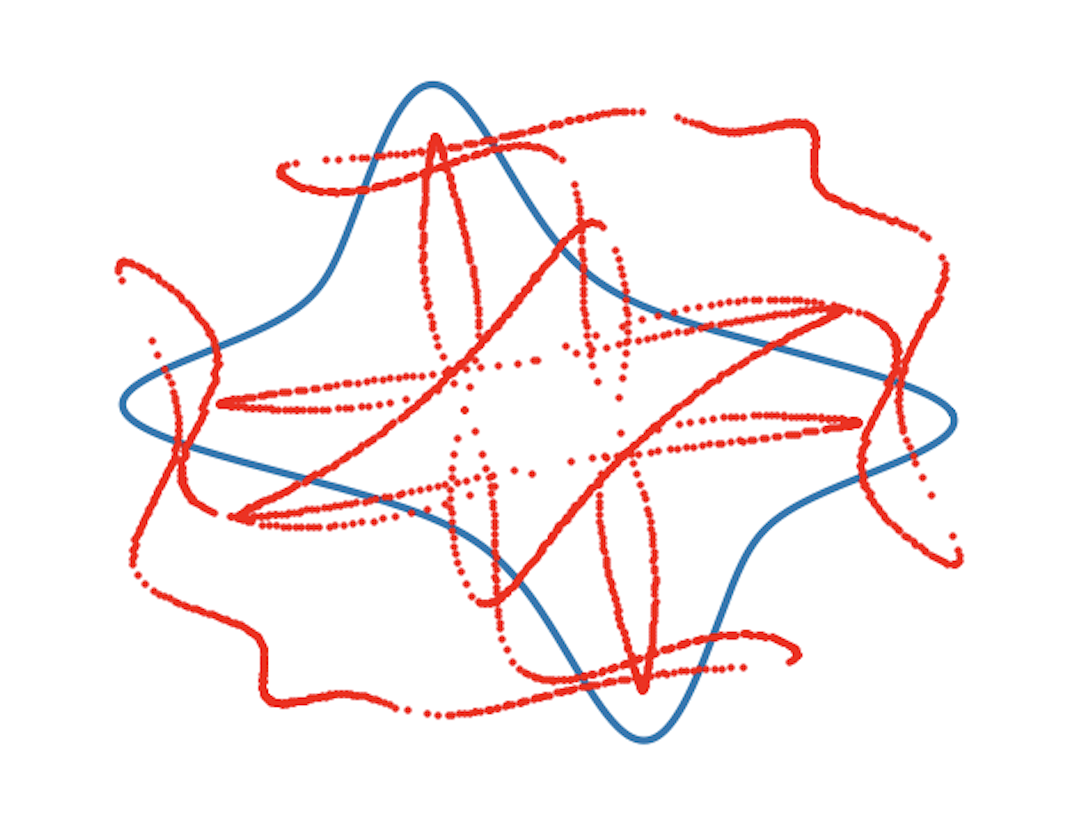}
\caption{In red a discretization of $\mathcal{SOSC}(\Gamma)$ where $\Gamma$ is the blue set \label{fig5}}
\end{center}
\end{figure}

Classical transversality arguments can be used to demonstrate\footnote{As a matter if fact, an appropriate Multijet Transversality Theorem (see {\it e.g.} \cite[Theorem 4.13 p. 57]{gg73}) allows to show that for a generic smooth simple closed curve, the set of $(t,u,v) \in \S^1$ such that $(\gamma(t),\gamma(u), \gamma(v)) \in \mathcal{SC}$ is a compact smooth submanifold of dimension $1$ of $(\S^1)^3$ and that the set of opposite square corners of the curve is the image of that set by a smooth immersion. } that the set of opposite square corners of a generic smooth simple closed curve is a (non necessarily connected\footnote{For example, we can show that if we consider the ellipse $\Gamma$ of Figure \ref{fig4} and replace for $\epsilon>0$ small the (short) piece of $\Gamma$  joining $A_{\epsilon}:=(1-\epsilon,2 \sqrt{2 \epsilon -\epsilon^2})$ to $B_{\epsilon}:=(1-\epsilon,2 \sqrt{2 \epsilon +\epsilon^2})$  by a (simple non-closed) curve from $A_{\epsilon}$ to $B_{\epsilon}$ contained in the ball centered at $(1,0)$ with radius $10\sqrt{\epsilon}$ and with a non-empty set of opposite square corners, then this small deformation generates a new connected component of $\mathcal{SOSC}(\Gamma)$ for $\epsilon$ small enough.}) compact smooth manifold of dimension $1$ immersed in the plane (see Figures \ref{fig4} and \ref{fig5}). Moreover, according to well-known results on the Square Peg Problem in the generic smooth case asserting that generic smooth (simple closed) curves inscribe an odd number of squares, one is inclined to think that the intersection of a set with its set of opposite square corners generically contains exactly a finite number of points which is an odd multiple of $4$. This conclusion follows from existing results which are based on purely topological methods that do not rely on the set of opposite square corners, and as we said before, this type of approach does not allow, {\it a priori}, any estimation on the sidelength of the squares in terms of the ``geometry"  of the curve. This is not the case of the method proposed in \cite{tao17}, where Tao proves a conservation lemma (see Lemma \ref{LEMTAO}) that allows to resolve the Square Peg Problem whenever the set of opposite square corners has a peculiar form. 

Before proceeding further, we mention that similar constructions of sets of opposite vertices have been used by Matschke in \cite{matschke20} to deal with the problem of quadrilaterals inscribed in convex curves. As in our paper, his results are based on a conservation lemma (see Lemma \ref{LEMTAO}) by Karasev \cite{karasev13} and Tao \cite{tao17} and considerations in terms of area of some sets related to the set of opposite vertices. We refer the reader to  \cite{matschke20} for further details.\\

\begin{figure}[H]
\begin{center}
\includegraphics[width=4cm]{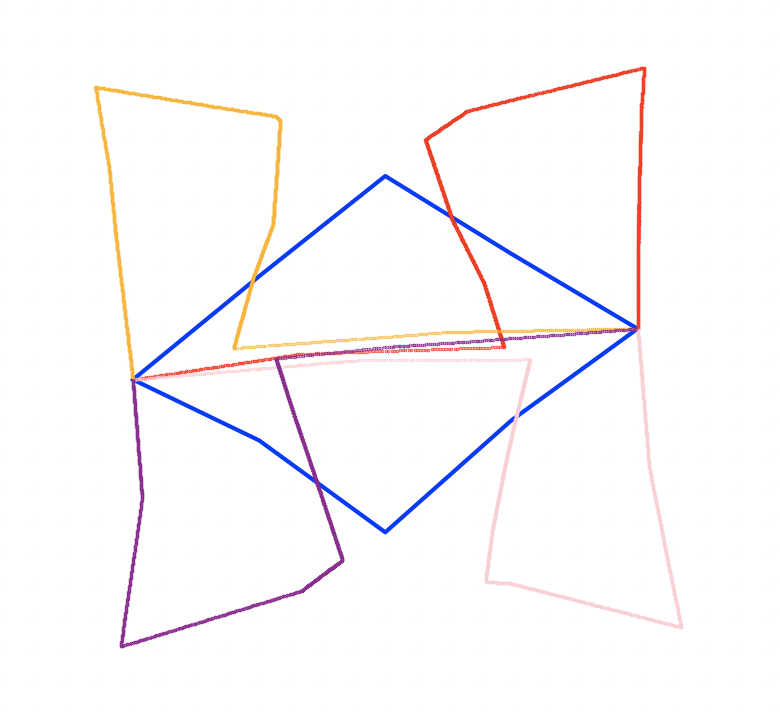}
\caption{The set  $\mathcal{SOSC}(\Gamma)$, with $\Gamma$ the set in blue, is the union of the four non-blue simple curves \label{fig6}}
\end{center}
\end{figure}

In \cite{tao17}, Tao considers simple closed curves given by the union of the graphs of two functions with ``small'' Lipschitz constants. Given an interval $I=[T_0,T_1]$, two functions $f,g:I \rightarrow \R$ such that 
\[
f(T_0)=g(T_0), \quad f(T_1)=g(T_1) \quad \mbox{and} \quad  f(t)<g(t) \, \, \forall t\in (T_0,T_1)
\]
and setting 
\[
\Gamma := \Gamma^f \cup \Gamma^g \quad \mbox{with} \quad \Gamma^f:=\mbox{Graph}_f (I), \,  \Gamma^g:=\mbox{Graph}_g (I),  
\]
he shows, roughly speaking\footnote{The result on $\mathcal{SOSC}(\Gamma)$ that we are stating here is not rigorously correct, we refer the reader to Tao's paper \cite{tao17} or to  Lemma \ref{LEM1} and Remarks \ref{REM1}, \ref{REM2} for a better understanding of the situation.}, that, if $f$ and $g$ are $(1-\epsilon)$-Lipschitz for some $\epsilon>0$, then the set $\mathcal{SOSC}(\Gamma)$ is the union of the four sets (see Figure \ref{fig6})
\[
\mathcal{S}_1  = \left\{ \mathcal{Q} (O,P,R) \, \vert \, (O,P,Q) \in \mathcal{SC}\cap \left( \Gamma^f \times \Gamma^f \times \Gamma^g\right) \right\},
\]
\[
\mathcal{S}_2  = \left\{ \mathcal{Q} (O,P,R) \, \vert \, (O,P,Q) \in \mathcal{SC}\cap  \left( \Gamma^f \times \Gamma^g \times \Gamma^f\right)  \right\},
\]
\[
\mathcal{S}_3  = \left\{ \mathcal{Q} (O,P,R) \, \vert \, (O,P,Q) \in \mathcal{SC}\cap  \left( \Gamma^g \times \Gamma^g \times \Gamma^f\right)  \right\},
\]
\[
\mathcal{S}_4  = \left\{ \mathcal{Q} (O,P,R) \, \vert \, (O,P,Q) \in \mathcal{SC}\cap  \left( \Gamma^g \times \Gamma^f \times \Gamma^g\right)  \right\},
\]

and moreover each of those sets is a Lipschitz simple curve joining the point $P_0:=(T_0,f(T_0))=(T_0,g(T_0))$ to the point $P_1:=(T_1,f(T_1))=(T_1,g(T_1))$. Then, Tao applies a conservation lemma (see Lemma \ref{LEMTAO})  to show that the (signed) area enclosed by $\Gamma^g$ and $\mathcal{S}_1$ has to be zero. Which implies that the two curves $\Gamma^g\setminus \{P_0,P_1\}$ and $\mathcal{S}_1$ must intersect and so proves\footnote{As we said, we sketch here a simplified version of Tao's proof, we refer the reader to \cite{tao17} for the complete proof.} the existence of a square inscribed in $\Gamma$. The idea of this present paper is simply to show that a ``quantification" of Tao's approach allows to obtain the following result:

\begin{theorem}\label{MainTHM}
There is a universal constant $C>0$ such that the following property holds: Let $I=[T_0,T_1]$ be interval and $f,g : I \rightarrow \R$ be $1$-Lipschitz functions such that $f(T_0)=g(T_0)$, $f(T_1)=g(T_1)$ and $f(t)<g(t)$ for all $t\in (T_0,T_1)$, then the set 
\[
\mbox{Graph}_f (I) \cup \mbox{Graph}_g (I)
\]
inscribes a square of sidelength at least
\[
C \cdot \max_{t\in I} \Bigl\{ g(t)-f(t) \Bigr\}.
\] 
\end{theorem}

The proof of Theorem \ref{MainTHM} is sketched in the next section and given with full detail in Section \ref{SECProof}. It provides a constant $C$ equal to $0.018$ which seems very far from being sharp since computer simulations\footnote{We wrote a computer program in Python to generate at random pairs $(f,g)$ of piecewise affine functions (over a dyadic partition of $[0,1]$) as in Theorem \ref{MainTHM} and compute the corresponding inscribed squares. We have not found any inscribed square with sidelength $<0.5 \cdot\max_{t\in I} \{ g(t)-f(t)\}$. By the way, we leave the reader to check the optimal constant has to be $\leq 0.5$. } suggest that the optimal constant of Theorem \ref{MainTHM} is probably $0.5$. Although moderately interesting, Theorem \ref{MainTHM} shows at least that Tao's approach might provide an efficient method to quantify the size of squares inscribed in simple closed curves in term of the geometry of the curve and thus might be certainly useful to settle the Square Peg Problem in the general case. But the road is long, there are a number of issues. As an example, the second part of the proof of Theorem \ref{MainTHM} relies heavily on the $1$-Lipschitzness assumption on the functions (see next section), can we weaken this assumption? More precisely, does a variant of Theorem \ref{MainTHM} (where $\max\{g-f\}$ can be replaced by another quantity depending on $f$ and $g$) holds true if we assume that $f$ ang $g$ are smooth (or by approximation only continuous), but not necessarily $1$-Lipschitz, and that the set $\mathcal{SOSC}(\Gamma)$ contains a simple Lipschitz curve connecting $(T_0,f(T_0)) = (T_0,g(T_0)) $ to $(T_1,f(T_1)) = (T_1,g(T_1)) $?\\

The paper is organized as follows: The main idea of the proof of Theorem \ref{MainTHM} is explained in Section \ref{SECIdea}, its complete proof is given in Section \ref{SECProof}, and technical lemmas are stated and proved in Section \ref{SECEst}. \\

{\bf Acknowledgements.} The author is indebted to B. Matschke and T. Tao for bringing several  references to his attention. 

\section{A rough idea of the proof of Theorem \ref{MainTHM}}\label{SECIdea}

\begin{figure}[H]
\centering
\begin{tikzpicture}[scale=1]
\draw[color=black] plot [smooth] coordinates {(-5,0.5)  (-4,0.25) (0,0) (4,0.5) (6,1)};
\draw (6,1) node[right]{$\gamma_1$};
\draw[color=blue] plot [smooth] coordinates { (-5,0) (-4+0.9*0.966,0.25-0.259*0.9) (0.939,0.342) (4+0.9,0.5) (6,0.25)};
\draw (6,0.25) node[right]{{\color{blue} $\gamma_2$}};
\draw[color=brown] plot [smooth] coordinates { (-5,1)  (-4+0.9*0.966+0.9*0.259,0.25-0.9*0.259+0.9*0.966)  (0.939-0.342,0.342+0.939) (4+0.9,0.5+0.9)  (6,1.8)};
\draw (6,1.8) node[right]{{\color{brown} $\gamma_3$}};
\draw[color=cyan] plot [smooth] coordinates { (-5,1.5)  (-4+0.9*0.259,0.25+0.9*0.966) (-0.342,0.939) (4,0.5+0.9)  (6,2.25)};
\draw (6,2.25) node[right]{{\color{cyan} $\gamma_4$}};
\fill[color=red] (0,0) circle (0.25mm);
\fill[color=red] (0.939,0.342) circle (0.25mm);
\fill[color=red] (0.939-0.342,0.342+0.939) circle (0.25mm);
\fill[color=red]  (-0.342,0.939) circle (0.25mm);
\draw[color=red] (0,0) -- (0.939,0.342);
\draw[color=red] (0.939,0.342) -- (0.939-0.342,0.342+0.939);
\draw[color=red] (0.939-0.342,0.342+0.939) -- (-0.342,0.939);
\draw[color=red] (0,0) -- (-0.342,0.939);
\fill[color=red] (-4,0.25) circle (0.25mm);
\fill[color=red] (-4+0.9*0.966,0.25-0.259*0.9) circle (0.25mm);
\fill[color=red] (-4+0.9*0.966+0.9*0.259,0.25-0.9*0.259+0.9*0.966) circle (0.25mm);
\fill[color=red] (-4+0.9*0.259,0.25+0.9*0.966) circle (0.25mm);
\draw[color=red] (-4,0.25) -- (-4+0.9*0.966,0.25-0.259*0.9);
\draw[color=red] (-4+0.9*0.966,0.25-0.259*0.9) -- (-4+0.9*0.966+0.9*0.259,0.25-0.9*0.259+0.9*0.966);
\draw[color=red]  (-4+0.9*0.966+0.9*0.259,0.25-0.9*0.259+0.9*0.966) -- (-4+0.9*0.259,0.25+0.9*0.966);
\draw[color=red] (-4,0.25) --  (-4+0.9*0.259,0.25+0.9*0.966);
\fill[color=red] (4,0.5) circle (0.25mm);
\fill[color=red] (4+0.9,0.5) circle (0.25mm);
\fill[color=red] (4+0.9,0.5+0.9) circle (0.25mm);
\fill[color=red]  (4,0.5+0.9) circle (0.25mm);
\draw[color=red] (4,0.5) -- (4+0.9,0.5);
\draw[color=red] (4+0.9,0.5) -- (4+0.9,0.5+0.9);
\draw[color=red] (4+0.9,0.5+0.9) --  (4,0.5+0.9);
\draw[color=red] (4,0.5) --  (4,0.5+0.9);
\end{tikzpicture}
\caption{For every $t$, $(\gamma_1(t)\gamma_2(t)\gamma_3(t)\gamma_4(t))$ is a square\label{fig7}}
\end{figure}
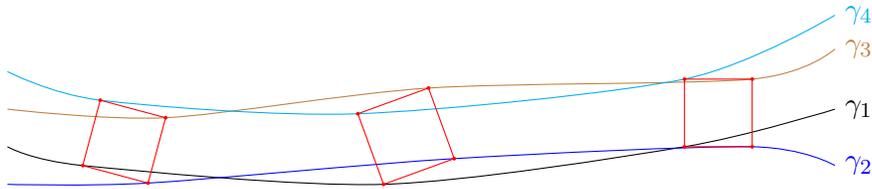

The proof of Theorem \ref{MainTHM} is based on two observations. The first one, due to Karasev \cite{karasev13} and Tao \cite{tao17}, is a result showing that some quantity is conserved along a quadruple of curves which traverse squares (see Figure \ref{fig7}). We refer the reader to \cite{tao17} for its proof and more details on the meaning of the integrals involved.

\vspace{0.2cm}
\vspace{0.2cm}

\begin{lemma}\label{LEMTAO}
Let $\gamma_1, \gamma_2, \gamma_3, \gamma_4 : [t_0,t_1] \rightarrow \R^2$ be rectifiables curves and $x,y,a,b: [t_0,t_1] \rightarrow \R$ be continuous functions such that 
\begin{eqnarray}\label{travsquare}
\begin{array}{rcl}
\gamma_1 (t) & = & \left( x(t), y(t) \right) \\
\gamma_2 (t) & = &  \left( x(t) + a(t), y(t) +b(t)\right)\\
\gamma_3 (t) & = &  \left( x(t) + a(t) -b(t), y(t) + a(t)+b(t)\right)\\
\gamma_4 (t) & = &  \left( x(t) - b(t), y(t) + a(t) \right)
\end{array}
\qquad \forall t \in [t_0,t_1].
\end{eqnarray}
Then we have the identity
\begin{eqnarray}
\int_{\gamma_1} y\, dx - \int_{\gamma_2} y\, dx + \int_{\gamma_3} y\, dx - \int_{\gamma_4} y\, dx = \frac{a(t_1)^2-b(t_1)^2}{2} - \frac{a(t_0)^2-b(t_0)^2}{2}. 
\end{eqnarray}
\end{lemma}

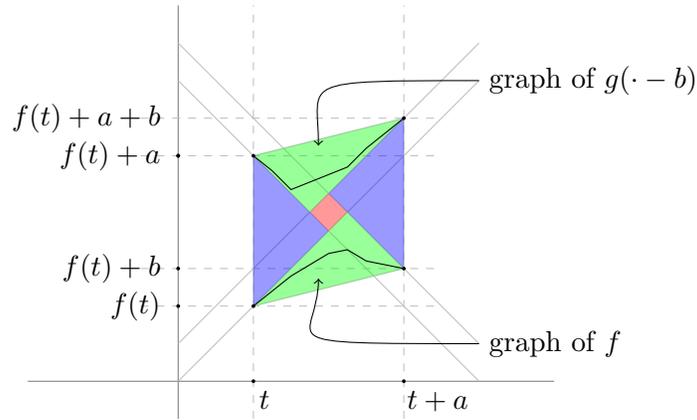
\begin{figure}[H]
\centering
\begin{tikzpicture}[scale=1]
\draw[color=gray] (0,-0.5) -- (0,5);
\draw[color=gray] (-2,0) -- (5,0);
\draw[color=lightgray, dashed] (1,-0.5) -- (1,5);
\fill[color=black] (1,0) circle (0.25mm);
\draw (1.15,0) node[below]{$t$};
\draw[color=lightgray, dashed] (3,-0.5) -- (3,5);
\fill[color=black] (3,0) circle (0.25mm);
\draw (3.45,0) node[below]{$t+a$};
\draw[color=lightgray, dashed] (-0.5,1) -- (3.5,1);
\fill[color=black] (0,1) circle (0.25mm);
\draw (-0.1,1) node[left]{$f(t)$};
\draw[color=lightgray, dashed] (-0.5,1.5) -- (3.5,1.5);
\fill[color=black] (0,1.5) circle (0.25mm);
\draw (-0.1,1.5) node[left]{$f(t)+b$};
\draw[color=lightgray, dashed] (-0.5,3) -- (3.5,3);
\fill[color=black] (0,3) circle (0.25mm);
\draw (-0.1,3) node[left]{$f(t)+a$};
\draw[color=lightgray, dashed] (-0.5,3.5) -- (3.5,3.5);
\fill[color=black] (0,3) circle (0.25mm);
\draw (-0.1,3.5) node[left]{$f(t)+a+b$};
\draw[color=lightgray] (1,1) -- (3,1.5) -- (3,3.5) -- (1,3) -- cycle; 
\draw[color=lightgray] (0,0.5) -- (4,4.5);
\draw[color=lightgray] (0,4) -- (4,0);
\draw[color=lightgray] (0,0) -- (4,4);
\draw[color=lightgray] (0,4.5) -- (4,0.5);
\fill[color=green, opacity=0.4] (1,3) -- (3,3.5) -- (1.75,2.25);
\fill[color=green, opacity=0.4] (1,1) -- (3,1.5) -- (2.25,2.25);
\fill[color=red, opacity=0.4] (1.75,2.25) -- (2,2.5) -- (2.25,2.25) -- (2,2) -- cycle;
\fill[color=blue, opacity=0.4] (1,3) -- (1,1) -- (2,2); 
\fill[color=blue, opacity=0.4] (3,3.5) -- (3,1.5) -- (2,2.5); 
\fill[color=black] (1,1) circle (0.25mm);
\fill[color=black] (3,1.5) circle (0.25mm);
\fill[color=black] (3,3.5) circle (0.25mm);
\fill[color=black] (1,3) circle (0.25mm);
\draw[color=black] (1,1) -- (1.5,1.4) -- (2,1.7) -- (2.25, 1.75) -- (2.5,1.6) -- (3,1.5);
\draw[color=black] (1,3) -- (1.25,2.8) -- (1.5,2.55) -- (2.25, 2.85) -- (2.5,3.1) -- (3,3.5);
\node[anchor=east] at (2,3) (text) {};
\node[anchor=west] at (4,4) (description) {graph of $g(\cdot-b)$};
\draw[->] (description) .. controls ([xshift=-4cm] description) and ([yshift=1cm] text) .. (text);
\node[anchor=east] at (2,1.5) (text) {};
\node[anchor=west] at (4,0.5) (description) {graph of $f$};
\draw[->] (description) .. controls ([xshift=-4cm] description) and ([yshift=-1cm] text) .. (text);
\end{tikzpicture}
\caption{$\int_{t-b}^{t+a-b} g(s)\, ds - \int_{t}^{t+a} f(s)\, ds=  \int_{t}^{t+a} g(s-b)-f(s)\, ds$ is larger or equal to the sum of the blue area and the red area given respectively by $a^2/2$ and $b^2/2$\label{fig8}}
\end{figure}

The above result will be used to show that if our simple closed curve (made of the union of the graphs of two  $1$-Lipschitz functions) does inscribe only squares with small sidelength then some area enclosed by a piece of graph of $g$ together with a piece of the set of opposite square corners has to be small. The second observation is the following type of result which gives a lower bound for some integral if the union of the graphs of two $1$-Lipschitz functions inscribes a certain type of square, its proof is given in Figure \ref{fig8}.

\begin{lemma}\label{LEMObs2}
Let $f,g : [T_0,T_1] \rightarrow \R$ be $1$-Lipschitz functions such that $f(T_0)=g(T_0)$, $f(T_1)=g(T_1)$ and $f(t)<g(t)$ for all $t\in (T_0,T_1)$ and $t,a,b \in \R$ be such that 
\begin{eqnarray}
T_0 \leq t, t+a, t-b, t+a-b \leq T_1, \quad  a>0
\end{eqnarray}
and
\begin{eqnarray}
\begin{array}{rcl}
f (t+a) & = & f(t)+b \\
g(t+a-b) & = &  f(t) +a+b\\
g(t-b)& =& f(t)+a.
\end{array}
\end{eqnarray}
Then we have 
\begin{eqnarray}
\int_{t-b}^{t+a-b} g(s)\, ds - \int_{t}^{t+a} f(s)\, ds \geq \frac{a^2+b^2}{2}.
\end{eqnarray}
\end{lemma}

\vspace{0.2cm}
\vspace{0.2cm}

Then the proof of Theorem \ref{MainTHM} is divided in two parts:\\

\noindent First part: We note that it is sufficient to prove the result for fonctions which are $(1-\epsilon)$-Lipschitz for some $\epsilon>0$ small. Then, we fix $\epsilon>0$ and two $(1-\epsilon)$-Lipschitz functions $f$ and $g$ as in the assumption of Theorem \ref{MainTHM}, we extend them to the whole real line and consider the set of (possibly degenerate) opposite square corners of the form 
\[
Q_t = P_t^f +R_t^g - O_t^f,
\]
where $t\in \R$, $O_t^f=(t,f(t))$, $P_t^f \in \Gamma^f$, $R_t^g\in \Gamma^g$ and $R_t^g$ is the image of $P_t^f$ by the rotation of angle $\pi/2$ with center $O_t^f$. As shown by Tao \cite{tao17}, the function $t\in \R \mapsto Q_t$ is Lipschitz and injective, so its image is a simple rectifiable curve joining $(T_0,f(T_0))$ to $(T_1,f(T_1))$, and by construction, the curve $t \in \R \mapsto Q_t$ comes along with three Lipschitz curves $t \in \R\mapsto O_t^f \in \Gamma^f, t \in \R\mapsto P_t^f \in \Gamma^f$ and $t\in \R\mapsto R_t^g \in \Gamma^g$ such that $(O_t^fP_t^fQ_tR_t^g)$ is always a square (see Lemma \ref{LEM1}). Therefore, any intersection of $\left(Q_t\right)_{t\in \R}$ with $\Gamma^g$ gives rise to a (possibly degenerate) inscribed square. Furthermore, if we consider two consecutive intersections of $\left(Q_t\right)_{t\in \R}$ with $\Gamma^g$, say at times $t_0<t_1$ in $[T_0,T_1]$ then Tao's conservation lemma can be used to bound from above the (non-signed) area enclosed by the simple closed curve made of the concatenation of $Q_{[t_0,t_1]}$ and the piece $\Gamma^g$ joining $Q_{t_0}$ to $Q_{t_1}$ in terms of the sizelengths of the corresponding squares at $t_0$ and $t_1$ (see Lemma \ref{LEM4}). This result allows to show that the smaller are the squares at $t_0$ and $t_1$ the smaller is the area and the closer to $\Gamma^g$ has to be the curve $t \in \R \mapsto Q_t$ (see Lemma \ref{LEM6}).\\

\noindent Second part: From the first part, we need to figure out what happens when the curve $t \in \R \mapsto Q_t$ is close to $\Gamma^g$ in some sense and to see how to get a contradiction if the squares at times $t_0$ and $t_1$ are too small. To do this let us imagine, for sake of simplicity,  that $Q_t$ is so close to $\Gamma^g$ that it belongs indeed to $\Gamma^g$ for all $t\in [t_0,t_1]$ and that the sidelength of the squares at $t_0, t_1$ is equal to $0$. Then in this case, Tao's Lemma \ref{LEMTAO} allows to show that we have (compare Lemma \ref{LEM5})
\begin{eqnarray*}
 \int_{t-b_{t}}^{t+a_t-b_t} g(s) \, ds -  \int_{t}^{t+a_t} f(s) \, ds=     \frac{a_{t}^2-b_{t}^2}{2} \qquad t\in [t_0,t_1],
\end{eqnarray*}
where we suppose that $O_t^f, P_t^f, Q_t$ and $R_t^g$ satisfy
\[
\begin{array}{rcl}
O_t^f & = & (t,f(t))\\
P_t^f & = &\left(t+a_t,f(t+a_t)\right) =  \left(t+a_t,f(t)+b_t\right)   \\
 Q_t &= & \left( t+a_t-b_t,f(t+a_t-b_t)\right) = \left( t+a_t-b_t,f(t)+a_t+b_t\right)\\
  R_t^g & = & \left(t-b_t,g(t-b_t)\right) = \left(t-b_t,f(t)+a_t\right)
  \end{array}
  \qquad \forall t \in [t_0,t_1].
\] 
Furthermore, Lemma \ref{LEMObs2} implies that we also have 
\begin{eqnarray*}
 \int_{t-b_{t}}^{t+a_t-b_t} g(s) \, ds -  \int_{t}^{t+a_t} f(s) \, ds \geq      \frac{a_{t}^2+b_{t}^2}{2} \qquad t\in [t_0,t_1].
\end{eqnarray*}
Then, we conclude that $b_t=0$ for all $t\in [t_0,t_1]$ and as a consequence that 
\[
g(t)=f(t)+a_t, \quad g(t+a_t)=f(t)+a_t \quad \mbox{and} \quad   \int_{t}^{t+a_t} g(s) -f(s) \, ds =     \frac{a_{t}^2}{2},
\]
for all $t\in [t_0,t_1]$. This type of property prevents the function $g-f$ to admit a maximum at some $T\in (t_0,t_1)$ such that $T\in (t,t+a_t) \in [t_0,t_1]$ for some $t\in [t_0,t_1]$, and yields a contradiction. \\

The proof of Theorem \ref{MainTHM} consists in quantifying all the above arguments to make the sketch of proof correct. 

\section{Proof of Theorem \ref{MainTHM}}\label{SECProof}

It is sufficient to prove the result for functions $f$ and $g$ which are $(1-\epsilon)$-Lipschitz for some $\epsilon>0$. As a matter of fact, if Theorem \ref{MainTHM} holds true in this case, then given $f,g:I \rightarrow \R$ we can define for every $\epsilon>0$ small $f_{\epsilon}, g_{\epsilon}:I \rightarrow \R$ by $f_{\epsilon}:=(1-\epsilon)f$ and $g_{\epsilon}:=(1-\epsilon)g$, apply the result and pass to the limit as $\epsilon \downarrow 0$ to obtain the required inscribed square as the limit of a sequence of squares whose sidelengths are bounded from below by $C\cdot \max_{t\in I}\{ g_{\epsilon}(t)-f_{\epsilon}(t) \}$ which tends to  $C\cdot \max_{t\in I}\{ g(t)-f(t) \}$. So from now on, we assume that we are given two functions $f, g: I \rightarrow \R$ which are $(1-\epsilon)$-Lipschitz for some $\epsilon>0$ and set
\[
M := \max_{t\in I} \Bigl\{ g(t)-f(t) \Bigr\} >0.
\]
Moreover, as in  \cite{tao17}, we extend $f$ and $g$ to the whole real line $\R$ by setting 
\[
f(t)=g(t)=f(T_0)=g(T_0) \, \forall t\leq T_0 \quad \mbox{and} \quad f(t)=g(t)=f(T_1)=g(T_1) \, \forall t\geq T_1,
\]
and denote respectively the graphs of $f$ and $g$ over $\R$ by $\Gamma^f$ and $\Gamma^g$. Then, for every $t\in \R$, we set $O_t^f=(t,f(t))$ and denote by $\mbox{Rot}^f_t:\R^2 \rightarrow \R^2$ the rotation of angle $\pi/2$ with center $O_t^f$. 

\begin{lemma}\label{LEM1}
The following properties hold:
\begin{itemize}
\item[(i)] For every $t\in \R$,  the set $\mbox{Rot}^f_t \left(\Gamma^f\right) \cap \Gamma^g$ is  a singleton equal to $\{R_t^g\}$ with $R_t^g= \mbox{Rot}^f_t (P_t^f) \in \Gamma^g$ and $P_t^f=(u_t,f(u_t)) \in \Gamma^f$, where $u_t \in [t+(g(t)-f(t))/2,+\infty)$ is the unique solution $u\in \R$ of the equation 
\[
g ( t+f(t)-f(u)) - f(t)-u+t =0.
\] 
\item[(ii)] If $t\notin (T_0,T_1)$, then  $\mbox{Rot}^f_t \left(\Gamma^f\right) \cap \Gamma^g = \{O_t^f\}= \{(t,g(t)\}$.
\item[(iii)] The functions $t\in \R \mapsto P_t^f$, $t\in \R \mapsto R_t^f$ and $Q:\R \rightarrow \R^2$ defined by
\[
Q(t)=Q_t := P_t^f + R_t^g-O_t^f \qquad \forall t \in \R
\]
are Lipschitz.
\item[(iv)] The Lipschitz functions $a, b:\R \rightarrow \R$ defined by $a_t := u_t-t$ and $b_t:=f(u_t)-f(t)$ for all $t\in \R$ satisfy
\[
\begin{array}{rcl}
P_t^f & = & \left(t+a_t,f(t)+b_t\right) = \left(t+a_t,f(t+a_t)\right)  \\
 Q_t &= & \left( t+a_t-b_t,f(t)+a_t+b_t\right)\\
  R_t^g & = & \left(t-b_t,f(t)+a_t\right) =  \left(t-b_t,g(t-b_t)\right) ,
  \end{array}
\] 
and
\[
a_t \geq \frac{g(t)-f(t)}{2}, \quad |b_t|\leq a_t
\]
for all $t\in \R$.
\item[(v)] The function $t\in \R\mapsto Q_t$ is injective.
\end{itemize}
\end{lemma}

\begin{remark}\label{REM1}
The triple $(O_t^f,P_t^f,R_t^g)$ is not always in  $\mathcal{SC}$ (the set of square corners), it is the case if and only if  $O_t^f\neq P_t^f \Leftrightarrow u_t \neq 0 \Leftrightarrow t \in (T_0,T_1)$. Moreover, we do not have necessarily $Q_t \in \mathcal{SOSC}(\mbox{Graph}_f (I) \cup \mbox{Graph}_g (I))$ for all $t\in (T_0,T_1)$ because we might have  that $P_t^f \in \Gamma^f \setminus \mbox{Graph}_f (I)$ (remember that $\Gamma^f$ denotes the graph of $f$ over $\R$) for some $t$ in $(T_0,T_1)$. 
\end{remark}
 
\begin{remark}\label{REM2}
If we denote for every $t\in \R$, by $\mbox{Rot}^{f,-}_t:\R^2 \rightarrow \R^2$ the rotation of angle $-\pi/2$ with center $O_t^f$, by $\mbox{Rot}^{g}_t:\R^2 \rightarrow \R^2$ the rotation of angle $\pi/2$ with center $O_t^g$ and by $\mbox{Rot}^{g,-}_t:\R^2 \rightarrow \R^2$ the rotation of angle $-\pi/2$ with center $O_t^g$, then Lemma \ref{LEM1} implies by symmetry (we can exchange the roles of $f$ anf $g$ and/or reverse time) that for every $t\in \R$, the sets $\mbox{Rot}^{f,-}_t \left(\Gamma^f\right) \cap \Gamma^g$, $\mbox{Rot}^g_t \left(\Gamma^g\right) \cap \Gamma^f$ and $\mbox{Rot}^{g,-}_t \left(\Gamma^g\right) \cap \Gamma^f$ are singletons and the corresponding mappings $t \in \R\mapsto Q_t^2$, $t \in \R\mapsto Q_t^3$ and $t\in \R \mapsto Q_t^4$ are Lipschitz and injective.  Moreover, we can check easily that 
\[
 \mathcal{SOSC}(\mbox{Graph}_f (I) \cup \mbox{Graph}_g (I)) \subset \bigcup_{t\in [T_0,T_1]} \left\{Q_t, Q_t^2,Q_t^3, Q_t^4\right\}.
\]
\end{remark}
 
\begin{remark}
Lemma \ref{LEM1} requires $f$ and $g$ to be $(1-\epsilon)$-Lipschitz for some $\epsilon>0$. If $f$ and $g$ are only $1$-Lipschitz then one can show that for every $t\in \R$, the set $\mbox{Rot}^f_t \left(\Gamma^f\right) \cap \Gamma^g$ is either a singleton or a segment of slope $\pm1$. 
\end{remark}

\begin{proof}[Proof of Lemma \ref{LEM1}]
Let $t\in \R$ be fixed. The continuous function $\varphi_t:[t,+\infty) \rightarrow \R$ defined by 
\[
\varphi_t(u) := g ( t+f(t)-f(u)) - f(t)-u+t  \qquad \forall u \in \R,
\]
satisfies 
\[
\lim_{u\rightarrow +\infty}\varphi_t(u)=-\infty
\]
and, by $1$-lipschitzness of $f$ and $g$ (they are indeed $(1-\epsilon)$-Lipschitz) and the non-negativity of $g-f$, we have
\begin{eqnarray*}
\varphi_t\left(t+\frac{g(t)-f(t)}{2}\right) & = & g \left(t+f(t)-f\left( t+\frac{g(t)-f(t)}{2}\right)\right) - \frac{f(t)}{2} - \frac{g(t)}{2}\\
& \geq & g(t) - \left| f(t)-f\left( t+\frac{g(t)-f(t)}{2}\right) \right|  - \frac{f(t)}{2} - \frac{g(t)}{2}\\
& \geq &  g(t) - \left| \frac{g(t)-f(t)}{2} \right|  - \frac{f(t)}{2} - \frac{g(t)}{2}=0.
\end{eqnarray*}
Hence there is $u\geq t+\frac{g(t)-f(t)}{2}$ such that $\varphi_t(u) =0$, that is, such that
\[
 \mbox{Rot}^f_t (u,f(u))= (t,f(t)) + ( f(t)-f(u), u-t)=(t+f(t)-f(u),f(t)+u-t) \in \Gamma^g.
\]
This $u$ is unique because if there is another $u'\in \R$ verifying $\varphi_t(u') =0$, then we have ( by $(1-\epsilon)$-lipschitzness of $f$ and $g$)
\begin{multline*}
|u'-u|  = \left|  g ( t+f(t)-f(u')) -  g ( t+f(t)-f(u)) \right| \\
 \leq  (1-\epsilon)\left| f(u')-f(u)\right| \leq (1-\epsilon)^2 |u'-u|,
\end{multline*}
which shows that $u'=u$. Thus the proof of (i) is complete. We notice that if $t\notin (T_0,T_1)$, then $\varphi_t(t)= g(t)-f(t)=0$, which shows that $u=t$, that is $R_t^g=P_t^f=O_t^f$, corresponds to the unique point in  $\mbox{Rot}^f_t \left(\Gamma^f\right) \cap \Gamma^g$ and 
proves (ii). For every $t, t'$ in $\R$, we have 
\begin{eqnarray*}
& \quad & \left| u_t-u_{t'}\right| \\
& = & \left| g\left(t+f(t)-f(u_t)\right) -f(t)+t - g\left(t'+f(t')-f(u_{t'})\right) +f(t')-t'\right|\\
& \leq & \left| g\left(t+f(t)-f(u_t)\right) - g\left(t'+f(t')-f(u_{t'})\right) \right| + \left| f(t')-f(t)\right| + \left| t'-t\right|\\
& \leq & (1-\epsilon) \left| \left(t+f(t)-f(u_t)\right) - \left(t'+f(t')-f(u_{t'})\right) \right| + (1-\epsilon) \left|t'-t\right| + \left| t'-t\right|\\
& \leq & \left| \left(t+f(t)-f(u_t)\right) - \left(t'+f(t')-f(u_{t'})\right) \right| + 2 \left| t'-t\right|\\
& \leq & \left| f(u_t)  -f(u_{t'}) \right| + \left| f(t')-f(t)\right| + 3 \left| t'-t\right|\\
& \leq & (1-\epsilon) \left| u_{t'}-u_t\right| + (1-\epsilon) \left|t'-t\right| + 3 \left| t'-t\right| \leq (1-\epsilon) \left| u_{t'}-u_t\right| + 4 \left| t'-t\right|,
\end{eqnarray*}
which implies that $|u_{t'}-u_t| \leq (4/\epsilon) |t'-t|$. This shows that the function $t\in \R \mapsto u_t$ is Lipschitz and as a consequence that the functions defined in (iii) are Lipschitz. The first part of (iv) is a straightforward consequence of the definitions of $P_t^f, R_t^g, Q_t$ and $a_t, b_t$. Concerning the second part, $a_t\geq (g(t)-f(t))/2$ follows from $u_t\geq t+ (g(t)-f(t))/2$ and $|b_t|\leq a_t$ follows from the $1$-lipschitzness of $f$ (because $f(t+a_t)=f(t)+b_t$). To prove (v) we suppose for contradiction that there are $t\neq t'\in \R$ such that $Q_t=Q_{t'}$. Then, the point $Q=Q_t=Q_{t'}$ belongs to the two squares 
\[
(O_t^fP_t^fQR_t^g) \quad \mbox{and} \quad (O_{t'}^fP_{t'}^fQR_{t'}^g),
\]
which shows that $P_t^f\neq P_{t'}^f$, $R_t^g\neq R_{t'}^g$ (otherwise $O_t^f=O_{t'}^f$ which is impossible because $t\neq t'$) and that $R_t^g$ (resp. $R_{t'}^g$) is the image of $P_t^f$ (resp. of $P_{t'}^f$)  by the rotation of angle $-\pi/2$ about $Q$. Therefore, the lines $D=(P_t^fP_{t'}^f)$ and $D'= (R_t^gR_{t'}^g)$ are well-defined (they pass through different points) and $D'$ is the image of $D$ by the rotation of angle $-\pi/2$ about $Q$. But since $P_t^f$ and $P_{t'}^f$ belong to $\Gamma^f$ and $f$ is $(1-\epsilon)$-Lipschitz the angle between $D$ and the horizontal is strictly less than $\pi/4$ and consequently the angle of $D'$, its image by the rotation of angle $-\pi/2$ about $Q$, with the vertical is strictly less than $\pi/4$. This is a contradiction because  $R_t^g$ and $R_{t'}^g$ belong to $\Gamma^g$ and $g$ is $(1-\epsilon)$-Lipschitz.
 \end{proof}

By construction, for every $t\in \R$, the points $O_t^f,P_t^f,Q_t$ and $R_t^g$ form a square, and in addition the points $O_t^f$ and $P_t^f$ always belong to $\Gamma^f$ and $R_t^g$ always belongs to $\Gamma^g$. So,  we can apply Lemma \ref{LEMTAO} and write some integrals in terms of integrals of $f$ and $g$. 

\begin{lemma}\label{LEM2}
For every $t<t' \in \R$, 
\begin{multline*}
 \int_{t}^{t+a_{t}} f(s) \, ds -  \int_{t'}^{t'+a_{t'}} f(s) \, ds + \int_{Q([t,t'])} y \, dx - \int_{t-b_{t}}^{t'-b_{t'}} g(s) \, ds =    \frac{a_{t'}^2-b_{t'}^2}{2} -   \frac{a_{t}^2-b_{t}^2}{2}.
\end{multline*}
\end{lemma}

\begin{proof}[Proof of Lemma \ref{LEM2}]
Applying Tao's Lemma \ref{LEMTAO} with $\gamma_1, \gamma_2, \gamma_3, \gamma_4 : [t,t'] \rightarrow \R^2$ defined by (see Lemma \ref{LEM1})
\begin{eqnarray*}
\begin{array}{rcl}
\gamma_1 (s) & = & O_s^f = \left( s, f(s) \right) \\
\gamma_2 (s) & = & P_s^f =  \left( s + a_s, f(s) +b_s\right) = \left(s+a_s,f(s+a_s)\right)\\
\gamma_3 (s) & = &  Q_s  = \left( s + a_s -b_s, f(s) + a_s+b_s \right)\\
\gamma_4 (s) & = & R_s^g= \left( s - b_s, f(s) + a_s \right) = \left(s-b_s,g(s-b_s)\right),
\end{array}
\qquad \forall s \in [t,t'],
\end{eqnarray*}
we obtain
\begin{eqnarray}\label{May7EQ1}
\int_{\gamma_1} y\, dx - \int_{\gamma_2} y\, dx + \int_{\gamma_3} y\, dx - \int_{\gamma_4} y\, dx = \frac{a_{t'}^2-b_{t'}^2}{2} - \frac{a_t^2-b_t^2}{2}. 
\end{eqnarray}
Since $\gamma_1, \gamma_2$ and $\gamma_4$ are graphs, we have (see \cite[Example 3.3]{tao17})
\[
\int_{\gamma_1} y\, dx = \int_t^{t'} f(s)\, ds, \quad  \int_{\gamma_2} y\, dx = \int_{t+a_t}^{t'+a_{t'}} f(s)\, ds, \quad \int_{\gamma_4} y\, dx =\int_{t-b_t}^{t'-b_{t'}} g(s)\, ds.
 \]
 Then (\ref{May7EQ1}) gives 
 \begin{eqnarray*}
 \int_t^{t'} f(s)\, ds -  \int_{t+a_t}^{t'+a_{t'}} f(s)\, ds + \int_{Q([t,t'])} y\, dx - \int_{t-b_t}^{t'-b_{t'}} g(s)\, ds = \frac{a_{t'}^2-b_{t'}^2}{2} - \frac{a_t^2-b_t^2}{2}
\end{eqnarray*}
and we conclude by the equality 
\[
 \int_{t}^{t'} f(s)\, ds =  \int_t^{t+a_t} f(s)\, ds + \int_{t+a_t}^{t'+a_{t'}} f(s)\, ds + \int_{t'+a_{t'}}^{t'} f(s)\, ds.
 \]
\end{proof}

Let $T \in (T_0,T_1)$ such that $(g-f)(T)=M$ be fixed and $\rho \in (0,1/8)$ a constant to be chosen later. We define $t_0,t_1 \in [T_0,T_1]$ by 
\[
t_0 = \max \Bigl\{ t  \in [T_0, T] \, \vert \, Q_t \in \Gamma^g \mbox{ and }  a_t\leq \rho M  \Bigr\} 
\]
and
\[
t_1 = \min \Bigl\{ t  \in [T,T_1] \, \vert \, Q_t \in \Gamma^g \mbox{ and }  a_t\leq \rho M  \Bigr\}, 
\]
which are well-defined because $a_{T_0}=0$ and $a_{T_1}=0$ by Lemma \ref{LEM1} (iv). Then we set
\[
\tau_0:=  t_0+a_{t_0}-b_{t_0} \quad \mbox{and} \quad \tau_1:=t_1+a_{t_1}-b_{t_1}
\]
and note that  by construction, the following result holds:

\begin{lemma}\label{LEM3}
We have 
\begin{eqnarray}\label{EQJune1_1}
t_0 \leq \tau_0 < T <  T+\frac{3M}{8} < t_1 \leq \tau_1
\end{eqnarray}
and
\begin{eqnarray}\label{EQJune1_3}
Q_{t_0} = \left( \tau_0,g(\tau_0)\right), \quad Q_{t_1} = \left( \tau_1,g(\tau_1)\right) \quad \mbox{and} \quad Q_t \notin \Gamma^g \qquad \forall t \in (t_0,t_1).
\end{eqnarray}
Furthermore, if we denote  by $\Omega \subset \R^2$ the bounded open set enclosed by the curve 
\[
\gamma : [0,t_1-t_0+\tau_1-\tau_0]  \longrightarrow  \R^2
\]
given by the concatenation 
of $\mathcal{C}:=Q([t_0,t_1])$ with the reversal of $\mbox{Graph}_g\left( [\tau_0, \tau_1]\right)$ (which is a simple closed curve thanks to the previous property), then there is  $\sigma\in \{-1,1\}$ such that the following property holds: For every $\tau\in (\tau_0,\tau_1)$, there are $\lambda_{\tau} >0$ and $t_{\tau}\in (t_0,t_1)$ such that 
\begin{eqnarray}\label{EQJune1_4}
(\tau,g(\tau)+\sigma \, \lambda) = Q_{t_{\tau}} \in \mathcal{C} \quad \mbox{and} \quad (\tau,g(\tau)+\sigma \, s) \in \Omega \qquad \forall s \in (0,\lambda_{\tau}).
\end{eqnarray}
Moreover, if $\sigma=1$ then the curve $\gamma$ is clockwise oriented and if $\sigma=-1$ it is anticlockwise oriented. 
\end{lemma}

\begin{proof}
By construction and the fact that both $a_{t_0}-b_{t_0}, a_{t_1}-b_{t_1}$ are nonnegative (see Lemma \ref{LEM1} (iv)), we already know that $t_0\leq T \leq t_1$, $t_0\leq \tau_0$ and $t_1\leq \tau_1$. Since $Q_{t_0}\in \Gamma^g$, we have, by Lemma \ref{LEM1} (iv), 
\[
g (\tau_0) = g\left(t_0+a_{t_0}-b_{t_0}\right) = f(t_0) + a_{t_0} + b_{t_0} = f\left(t_0+a_{t_0}\right) + a_{t_0},
\]
which gives (by $1$-lipschitzness of $f$)
\begin{eqnarray}\label{EQJune1_2}
\left| g(\tau_0) - f(\tau_0) \right| \leq \left| g(\tau_0) - f(t_0 +a_{t_0} ) \right| +  \left|  f(t_0 +a_{t_0} ) - f(\tau_0)\right| \leq 2 a_{t_0} \leq 2 \rho M.
\end{eqnarray}
Consequently, if $\tau_0\geq T$, then we have (remember that $|b_{t_0}|\leq a_{t_0}$)
\[
0 \leq \tau_0-T \leq \tau_0-t_0 = a_{t_0}-b_{t_0} \leq 2 a_{t_0} \leq 2 \rho M,
\]
which implies (by $2$-lipschitzness of $g-f$)
\[
\left(g(\tau_0) -f(\tau_0)\right) \geq \left(g(T) -f(T)\right) - 2\left| \tau_0-T \right| \geq M - 4 \rho M = (1-4\rho) \, M,
\]
which contradicts (\ref{EQJune1_2}) since $\rho<1/8$. Furthermore, since $g(t_1-b_1)=f(t_1)+a_{t_1}$ and $|b_{t_1}|\leq a_{t_1}$, we have (by $1$-lipschitzness of $f$ and $g$)
\begin{multline*}
M - 2 (t_1-T) = g(T)-f(T) - 2 (t_1-T)  \leq \\
g(t_1)-f(t_1) \\
\leq  \left| g(t_1)  - g(t_1-b_{t_1})\right| + \left| g(t_1-b_{t_1})-f(t_1) \right| \leq  \ 2a_{t_1} \leq 2\rho M
\end{multline*}
which implies that $t_1-T\geq M(1-2\rho)/2 > 3M/8$ since $\rho\in (0,1/8)$. So, the proof of (\ref{EQJune1_1}) is complete. The property (\ref{EQJune1_3}) is a direct consequence of the construction of $t_0$ and $t_1$. Let us now prove the second part of the statement. The concatenation of $\mathcal{C}:=Q([t_0,t_1])$ with the reversal of $\mbox{Graph}_g\left( [\tau_0, \tau_1]\right)$ is the curve 
\[
\gamma : [0,t_1-t_0+\tau_1-\tau_0]  \longrightarrow  \R^2
\]
defined by 
\[
\gamma(s) := \left\{ \begin{array}{cl}
Q_{t_0+s} & \mbox{ if } s \in [0,t_1-t_0]\\
\left(\tau_1 - (s-t_1+t_0), g\left(  \tau_1 - (s-t_1+t_0)\right) \right) & \mbox{ if } s \in [t_1-t_0,t_1-t_0+\tau_1-\tau_0],
\end{array}
\right.
\]
for all $s\in [0,t_1-t_0+\tau_1-\tau_0]$. We check easily that $\gamma$ is Lipschitz and closed as the concatenation of two Lipschitz curves with the same endpoints (we have, by (\ref{EQJune1_3}), $\gamma(0)=Q_{t_0}=\left( \tau_0,g(\tau_0)\right)$ and $\gamma(t_1-t_0) = Q_{t_1}=\left( \tau_1,g(\tau_1)\right)$) and that it is simple because $t\in \R \mapsto Q_t$ is injective and $Q_t \notin \Gamma^g$ for all $t\in [t_0,t_1]$ (by (\ref{EQJune1_3})). Then, by the Jordan curve Theorem, the image of $\gamma$, $\mathcal{I}:= \gamma ([0,t_1-t_0+\tau_1-\tau_0])$,  divides the plane $\R^2$ into two connected components $\Omega$ and $O$ where $\Omega$ is the bounded open set enclosed by $\gamma$ and $O$ is the complement of $\overline{\Omega}$. For every $\tau\in (\tau_0,\tau_1)$, the point $(\tau,g(\tau))$ belongs to the image of $\gamma$ but not to $\gamma([0,t_1-t_0])$. So, since $g$ is $1$-Lipschitz there is $h_{\tau}>0$ such that the vertical segment centered at $(\tau,g(\tau))$ of length $2h_{\tau}$ intersects  $\mathcal{I}$ only at $(\tau,g(\tau))$ (note that the vertical line through $(\tau,g(\tau))$ intersects $\mathcal{I}$ at $(\tau,g(\tau))$ and at at least another point of $\mathcal{C}$ either at h . By connectedness of $\Omega$ and $O$ (and the fact that they are separated by $\mathcal{I}$), we can show that there is $\sigma\in \{-1,+1\}$ such that for every $\tau \in (\tau_0,\tau_1)$, $(\tau,g(\tau))+\sigma (0,h_{\tau})$ belongs to $\Omega$ and we can check that $\gamma$ is clockwise oriented if $\sigma=1$ and anticlockwise oriented if $\sigma=-1$. If for every $\tau \in (\tau_0,\tau_1)$, we define $\lambda_{\tau}>0$ by 
\[
\lambda_{\tau} := \min \Bigl\{ h>0 \, \vert \, (\tau,g(\tau)) + \sigma (0,h) \notin \Omega \Bigr\},
\]
then, since $\Omega$ is bounded $\lambda_{\tau}$ is well-defined, since $\lambda_{\tau}\geq h_{\tau}$ $\lambda_{\tau}$ is positive, and by construction $(\tau,g(\tau)) + \sigma (0,s) \in \Omega$ for $s\in (0,\lambda_{\tau})$ and $(\tau,g(\tau)) + \sigma (0,\lambda_{\tau})\in \overline{\Omega} \setminus \Gamma^g = \mathcal{C}$ so that there is a unique $t_{\tau}\in [t_0,t_1]$ (by  Lemma \ref{LEM1} (v)) such that  $(\tau,g(\tau)) + \sigma (0,\lambda_{\tau}) = Q_{t_{\tau}}$. This completes the proof of the Lemma. 
\end{proof}

By Lemma \ref{LEM3}, we know that for every $\tau \in (\tau_0,\tau_1)$ the concatenation of the curve $\mathcal{C}_{\tau}:=Q([t_0,t_{\tau}])$ with the vertical segment joining $Q_{t_{\tau}}$ to $(\tau,g(\tau))$ and the reversal of   $\mbox{Graph}_g\left( [\tau_0, \tau]\right)$ is a simple closed curve, we denote by $\Omega_{\tau}$ the open set enclosed by that curve and we set 
\[
A_{\tau} =  \int_{Q([t_0,t_{\tau}])} y \, dx - \int_{\tau_0}^{\tau} g(s)\, ds.
\]

\begin{figure}[H]
\centering
\begin{tikzpicture}[scale=2]
\foreach \p in {0,...,100}
{
	\draw[color=blue] ( {(-4*(0.01*\p)*(0.01*\p)+4*(0.01*\p)+1)*((0.01*\p)-0.5)+0.5},  { (-4*(0.01*\p)*(0.01*\p)+4*(0.01*\p)+1)*(sin(7.2*\p)*sin(7.2*\p)*0.7+0.005*\p*(1-0.01*\p) - 0.01*0.01*0.01*0.01*\p*\p*\p*\p*(1-0.01*\p)*(1-0.01*\p)  )} ) -- ( {(-4*(0.01*(\p+1))*(0.01*(\p+1))+4*(0.01*(\p+1))+1)*((0.01*(\p+1))-0.5)+0.5},  { (-4*(0.01*(\p+1))*(0.01*(\p+1))+4*(0.01*(\p+1))+1)*(sin(7.2*(\p+1))*sin(7.2*(\p+1))*0.7+0.005*(\p+1)*(1-0.01*(\p+1)) - 0.01*0.01*0.01*0.01*(\p+1)*(\p+1)*(\p+1)*(\p+1)*(1-0.01*(\p+1))*(1-0.01*(\p+1))  )} ) ;
}
\fill[opacity=0.2,blue] (0,0) \foreach \p in {1,...,54}{ --   ( {(-4*(0.01*\p)*(0.01*\p)+4*(0.01*\p)+1)*((0.01*\p)-0.5)+0.5},  { (-4*(0.01*\p)*(0.01*\p)+4*(0.01*\p)+1)*(sin(7.2*\p)*sin(7.2*\p)*0.7+0.005*\p*(1-0.01*\p) - 0.01*0.01*0.01*0.01*\p*\p*\p*\p*(1-0.01*\p)*(1-0.01*\p)  )} ) } -- (0.57,0.071) -- (0.4,0.1) -- cycle;
\draw[color=black] (-0.4,-0.2) -- (0,0) -- (0.4,0.1) -- (1,0) -- (1.4,0.1); 
\node[right] at (1.5,0) {graph of $g$};
\node[right] at (1.5,1) {\color{blue} curve $\mathcal{C}$};
\draw[color=gray] (-1,-0.5) -- (2.5,-0.5);
\draw (0,-0.5) node[below]{$\tau_0$};
\draw (0.57,-0.5) node[below]{$\tau$};
\draw (1,-0.5) node[below]{$\tau_1$};
\draw[color=lightgray, dashed] (0,-0.5) -- (0,2);
\draw[color=lightgray, dashed] (0.57,-0.5) -- (0.57,2);
\draw[color=lightgray, dashed] (1,-0.5) -- (1,2);
\draw[line width=0.05mm, blue] (0.57,0.5186) -- (0.57,0.071);
\fill[color=black] (0,-0.5) circle (0.125mm);
\fill[color=black] (0.57,-0.5) circle (0.125mm);
\fill[color=black] (1,-0.5) circle (0.125mm);
\fill[color=blue] (0.57,0.5186) circle (0.125mm);
\node[right] at (0.55,0.518) {\color{blue} $Q_{t_{\tau}}$};
\node[anchor=east] at (0.25,0.3) (text) {};
\node[anchor=east] at (-1,0.3) (description) {\color{blue} $\Omega_{\tau}$};
\draw[dashed,blue,->] (description) .. controls ([xshift=1cm] description) and ([yshift=0cm] text) .. (text);
\end{tikzpicture}
\caption{The set $\Omega_{\tau}$ is contained in $\Omega$, the set enclosed by the concatenation of $\mathcal{C}$ with $\mbox{Graph}_g([\tau_0, \tau_1])$ \label{fig9}}
\end{figure}
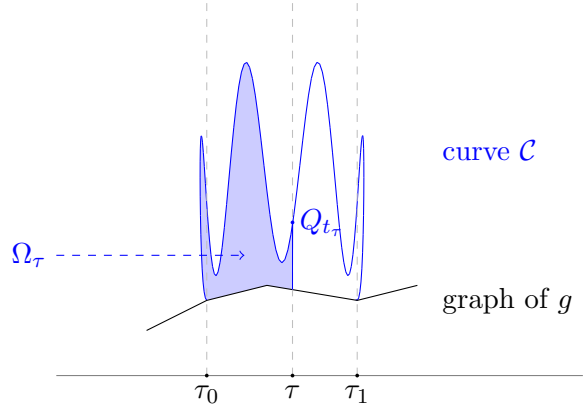

The following result follows from Stoke's formula, Lemma \ref{LEM2} and Lemma  \ref{LEMEst1} whose proof can be found in Section \ref{SECEst}.
 
\begin{lemma}\label{LEM4}
We have 
\begin{eqnarray}
\mathcal{L}^2(\Omega_{\tau}) = \left| A_{\tau} \right| \leq \mathcal{L}^2(\Omega) \leq \rho^2 M^2 \qquad \forall \tau \in (\tau_0,\tau_1).
\end{eqnarray}
Furthermore, if $\sigma=1$, then we have $0\leq A_{\tau}\leq A_{\tau_1}$ for all $\tau \in [\tau_0,\tau_1]$ and if  $\sigma=-1$, then we have $A_{\tau_1} \leq A_{\tau}\leq 0$ for all $\tau \in [\tau_0,\tau_1]$.
\end{lemma}
 
\begin{proof}[Proof of Lemma \ref{LEM4}]
By Lemma \ref{LEM3} and Stoke's formula, we have (see \cite[Lemma 3.4]{tao17})
\begin{eqnarray}\label{June3_1}
\mathcal{L}^2(\Omega) = \sigma \left( \int_{Q([t_0,t_1])} y \, dx - \int_{\tau_0}^{\tau_1} g(s)\, ds \right) = \sigma A_{\tau_1}
\end{eqnarray}
and since for every $\tau\in (\tau_0,\tau_1)$ the concatenation of $\mathcal{C}_{\tau}$ with the vertical segment joining $Q_{t_{\tau}}$ to $(\tau,g(\tau))$ and the reversal of $\mbox{Graph}_g\left( [\tau_0, \tau]\right)$ is a simple closed curve with the same orientation of $\gamma$, we also have 
\begin{eqnarray}\label{June3_2}
\mathcal{L}^2(\Omega_{\tau}) = \sigma A_{\tau} \qquad \forall \tau\in (\tau_0,\tau_1).
\end{eqnarray}
By construction, $\Omega_{\tau}$ is contained in $\Omega$ for all $\tau\in (\tau_0,\tau_1)$, so we have $\mathcal{L}^2(\Omega_{\tau}) = \left| A_{\tau} \right| \leq \mathcal{L}^2(\Omega)$ for every $\tau\in (\tau_0,\tau_1)$ and the second part of the lemma follows easily from (\ref{June3_1})-(\ref{June3_2}). It remains to show that $\mathcal{L}^2(\Omega)\leq \rho^2M^2$. We note that by Lemma \ref{LEM2}, we have
 \begin{eqnarray*}
& \quad & \int_{Q([t_0,t_1])} y \, dx - \int_{t_0+a_{t_0}-b_{t_0}}^{t_1+a_{t_1}-b_{t_1}} g(s)\, ds \\
& = &   \int_{Q([t_0,t_1])} y \, dx - \int_{t_0-b_{t_0}}^{t_1-b_{t_1}} g(s)\, ds -  \int_{t_0+a_{t_0}-b_{t_0}}^{t_0-b_{t_0}} g(s)\, ds -  \int_{t_1-b_{t_1}}^{t_1+a_{t_1}-b_{t_1}} g(s)\, ds\\
& = &    \frac{a_{t_1}^2-b_{t_1}^2}{2} -   \frac{a_{t_0}^2-b_{t_0}^2}{2} - \int_{t_0}^{t_0+a_{t_0}} f(s) \, ds +  \int_{t_1}^{t_1+a_{t_1}} f(s) \, ds \\
& \quad & \qquad \qquad \qquad \qquad \qquad \qquad   -  \int_{t_0+a_{t_0}-b_{t_0}}^{t_0-b_{t_0}} g(s)\, ds -  \int_{t_1-b_{t_1}}^{t_1+a_{t_1}-b_{t_1}} g(s)\, ds\\
& = &   \frac{a_{t_1}^2-b_{t_1}^2}{2} -   \frac{a_{t_0}^2-b_{t_0}^2}{2} + \left[  \int_{t_0-b_{t_0}}^{t_0+a_{t_0}-b_{t_0}} g(s)\, ds - \int_{t_0}^{t_0+a_{t_0}} f(s) \, ds  \right] \\
& \quad & \qquad \qquad \qquad \qquad \qquad \qquad   - \left[ \int_{t_1-b_{t_1}}^{t_1+a_{t_1}-b_{t_1}} g(s)\, ds -\int_{t_1}^{t_1+a_{t_1}} f(s) \, ds\right].
 \end{eqnarray*}
 But Lemma \ref{LEMEst1} gives 
 \[
 \frac{a_{t_i}^2+b_{t_i}^2}{2} \leq  \int_{t_i-b_{t_i}}^{t_i+a_{t_i}-b_{t_i}} g(s)\, ds -\int_{t_i}^{t_i+a_{t_i}} f(s) \, ds \leq \frac{3a_{t_i}^2-b_{t_i}^2}{2} \qquad \forall i=0,1.
 \]
 Then, we conclude that
 \[ 
-\rho^2M^2 \leq -a_{t_1}^2 \leq-a_{t_1}^2 + b_{t_0}^2 \leq \int_{Q([t_0,t_1])} y \, dx - \int_{\tau_0}^{\tau_1} g(s)\, ds \leq a_{t_0}^2 - b_{t_1}^2 \leq a_{t_0}^2\leq \rho^2M^2,
 \]
 which completes the proof of the lemma.
\end{proof}

\begin{lemma}\label{LEM5}
We have for every $\tau \in [\tau_0,\tau_1]$
\[
\int_{ t_{\tau}-b_{t_{\tau}}}^{\tau} g(s) \, ds - \int_{t_{\tau}}^{t_{\tau}+a_{t_{\tau}} } f (s) \, ds \leq \frac{a_{t_{\tau}}^2-b_{t_{\tau}}^2}{2} +  \rho^2 M^2.
 \]
 \end{lemma}
 
\begin{proof}[Proof of Lemma \ref{LEM5}]
Let $\tau \in [\tau_0,\tau_1]$ be fixed. If $\sigma=1$, then Lemma \ref{LEM2} with $t=t_0$ and $t'=t_{\tau}$, the definition of $A_{\tau}$ and  $A_{\tau}\geq 0$ (by Lemma \ref{LEM4}) give
\begin{eqnarray*}
& \quad & \int_{ t_{\tau}-b_{t_{\tau}}}^{\tau} g(s) \, ds - \int_{t_{\tau}}^{t_{\tau}+a_{t_{\tau}} } f (s) \, ds \\
& = & \int_{ t_{\tau}-b_{t_{\tau}}}^{\tau} g(s) \, ds +  \frac{a_{t_{\tau}}^2-b_{t_{\tau}}^2}{2} -  \frac{a_{t_{0}}^2-b_{t_{0}}^2}{2} -  \int_{t_{0}}^{t_{0}+a_{t_{0}} } f (s) \, ds\\
& \quad & \qquad \qquad \qquad \qquad \qquad \qquad  -  \int_{Q([t_0,t_{\tau}])} y \, dx + \int_{t_0-b_{t_0}}^{t_{\tau}-b_{t_{\tau}}} g(s)\, ds \\
& = &  \frac{a_{t_{\tau}}^2-b_{t_{\tau}}^2}{2} -  \frac{a_{t_{0}}^2-b_{t_{0}}^2}{2}  +  \int_{ t_{0}-b_{t_{0}}}^{\tau_0} g(s) \, ds  -  \int_{t_{0}}^{t_{0}+a_{t_{0}} } f (s) \, ds  - A_{\tau} \\
& \leq &  \frac{a_{t_{\tau}}^2-b_{t_{\tau}}^2}{2} -  \frac{a_{t_{0}}^2-b_{t_{0}}^2}{2}  +  \int_{ t_{0}-b_{t_{0}}}^{\tau_0} g(s) \, ds  -  \int_{t_{0}}^{t_{0}+a_{t_{0}} } f (s) \, ds,
\end{eqnarray*}
which can be bounded from above, by (\ref{Est2}) of Lemma \ref{LEMEst1} with $\delta=0$, by 
\[
  \frac{a_{t_{\tau}}^2-b_{t_{\tau}}^2}{2} -  \frac{a_{t_{0}}^2-b_{t_{0}}^2}{2}  +  \frac{3a_{t_{0}}^2-b_{t_{0}}^2}{2} =  \frac{a_{t_{\tau}}^2-b_{t_{\tau}}^2}{2} + a_{t_0}^2 \leq  \frac{a_{t_{\tau}}^2-b_{t_{\tau}}^2}{2}  + \rho^2 M^2.
\]
If $\sigma=-1$, then Lemma \ref{LEM2} with $t=t_{\tau}$ and $t'=t_1$ gives
\begin{multline*}
\int_{ t_{\tau}-b_{t_{\tau}}}^{\tau} g(s) \, ds - \int_{t_{\tau}}^{t_{\tau}+a_{t_{\tau}} } f (s) \, ds = \\
  \frac{a_{t_{\tau}}^2-b_{t_{\tau}}^2}{2} -  \frac{a_{t_{1}}^2-b_{t_{1}}^2}{2}   -  \int_{t_{1}}^{t_{1}+a_{t_{1}} } f (s) \, ds  +  \int_{Q([t_{\tau},t_{1}])} y \, dx - \int_{\tau}^{t_1-b_{t_1} } g(s) \, ds
\end{multline*}
where we have 
\begin{multline*}
 \int_{Q([t_{\tau},t_{1}])} y \, dx - \int_{\tau}^{t_1-b_{t_1} } g(s) \, ds =  \int_{Q([t_0,t_{1}])} y \, dx -  \int_{Q([t_0,t_{\tau}])} y \, dx \\
 - \int_{\tau}^{\tau_0 } g(s) \, ds  - \int_{\tau_0}^{\tau_1 } g(s) \, ds  - \int_{\tau_1}^{t_1-b_{t_1} } g(s) \, ds= A_{\tau_1} - A_{\tau} + \int_{t_1-b_{t_1} }^{\tau_1} g(s) \, ds.
 \end{multline*}
Therefore, since $A_{\tau_1} \leq A_{\tau} \leq 0$ (by Lemma \ref{LEM4}), we infer that
\begin{multline*}
\int_{ t_{\tau}-b_{t_{\tau}}}^{\tau} g(s) \, ds - \int_{t_{\tau}}^{t_{\tau}+a_{t_{\tau}} } f (s) \, ds \\
\leq  \frac{a_{t_{\tau}}^2-b_{t_{\tau}}^2}{2} -  \frac{a_{t_{1}}^2-b_{t_{1}}^2}{2}   -  \int_{t_{1}}^{t_{1}+a_{t_{1}} } f (s) \, ds +  \int_{t_1-b_{t_1} }^{\tau_1} g(s) \, ds,
\end{multline*}
which can be bounded from above, by (\ref{Est2}) of Lemma \ref{LEMEst1} with $\delta=0$, by 
\[
  \frac{a_{t_{\tau}}^2-b_{t_{\tau}}^2}{2} -  \frac{a_{t_{1}}^2-b_{t_{1}}^2}{2}  +  \frac{3a_{t_{1}}^2-b_{t_{1}}^2}{2} =  \frac{a_{t_{\tau}}^2-b_{t_{\tau}}^2}{2} + a_{t_1}^2 \leq  \frac{a_{t_{\tau}}^2-b_{t_{\tau}}^2}{2}  + \rho^2 M^2.
\]
\end{proof}

For every $\mu>0$, we set
\[
\Lambda_{\mu} := \Bigl\{ \tau \in [\tau_0,\tau_1] \, \vert \, \lambda_{\tau} \geq \mu\Bigr\}.
\]
The Lebesgue measure of this set is controlled by the area of $\Omega$, we have:

\begin{lemma}\label{LEM6}
For every $\mu>0$, $\mathcal{L}^1(\Lambda_{\mu}) \leq \mathcal{L}^2 (\Omega)/\mu$.
\end{lemma}

\begin{proof}[Proof of Lemma \ref{LEM6}]
By Fubini's Theorem, we have 
\[
\mathcal{L}^2(\Omega) = \int_{\R} \mathcal{H}^1 (\Omega \cap V_{\tau}) \, d\tau,
\]
where $V_{\tau}$ denotes the vertical line of abscissa $\tau$. If $\tau \in [\tau_0,\tau_1]$ is such that $\lambda_{\tau} \geq \mu$, then by (\ref{EQJune1_4}) of Lemma \ref{LEM3}, the $1$-dimensional set $\Omega \cap V_{\tau}$ contains at least a vertical segment of length $\mu$, so that  $ \mathcal{H}^1 (\Omega \cap V_{\tau}) \geq \mu$. As a consequence, we have that 
\[
\mathcal{L}^2(\Omega) \geq \int_{\Lambda_{\mu}} \mu \, d\tau = \mu \, \mathcal{L}^1\left( \Lambda_{\mu} \right),
\]
which proves the result. 
\end{proof}

We are now ready to conclude the proof of Theorem \ref{MainTHM}. We consider a constant $B\geq 2$ to be fixed later and set 
\[
\nu:=\rho M \quad \mbox{and} \quad \mu:=2B \rho^2M.  
\]
Since 
\[
 \left[ T+ \frac{M}{4} - \frac{M}{4B},  T+ \frac{M}{4}+\frac{M}{4B} \right] \subset \left[ \tau_0,\tau_1\right] \quad (\mbox{by (\ref{EQJune1_1})})
\]
and
\[
\mathcal{L}^1 \left(  \left[ T+ \frac{M}{4} - \frac{M}{4B},  T+ \frac{M}{4}+\frac{M}{4B} \right] \right) =   \frac{M}{2B}= \frac{\nu^2}{\mu},
\]
there exists by Lemma \ref{LEM4} ($\mathcal{L}^2(\Omega) \leq \nu^2$) and Lemma \ref{LEM6}
\begin{eqnarray}\label{20May0}
\tau \in \left[ T+ \frac{M}{4} -\frac{M}{4B},  T+ \frac{M}{4} +\frac{M}{4B} \right]
\end{eqnarray}
such that
\begin{eqnarray}\label{20May5}
\lambda_{\tau} \in [0, \mu].
\end{eqnarray}
We set $\lambda := \lambda_{\tau}, t:=t_{\tau}, a:=a_{t_{\tau}}, b:= b_{t_{\tau}}$ and
\[
 E:=  \int_{ t-b}^{\tau} g(s) \, ds - \int_{t}^{t+a} f (s) \, ds.
 \]
 The contradiction will come from two inconsistent bounds for $E$, one from above given by Lemma \ref{LEM5} and one from below that will follow from Lemma \ref{LEMEst2}. On the one hand, Lemma \ref{LEM5} gives
 \begin{eqnarray}\label{June2_1}
 E \leq \frac{a^2-b^2}{2} +  \nu^2.
 \end{eqnarray}
 Now, in order to apply Lemma \ref{LEMEst2}, we need to show that $\tau-a<T$. Let us do it. By Lemma \ref{LEM1} (iv), we have
\begin{eqnarray}\label{20May1}
\begin{array}{rcl}
\tau & = & t+a-b\\
f(t+a) & = & f(t)+b\\
g(t+a-b) & = & f(t)+a+b-\sigma \lambda\\
g(t-b) & = & f(t)+a,
  \end{array}
\end{eqnarray}
which can be used to show that (by using the $1$-lipschitzness of $f$)
\begin{eqnarray}\label{20May3}
0 < a & = & g(t+a-b) - f(t+a) + \sigma \lambda \nonumber\\
& \leq & g(\tau) - f(\tau)  + \left| f(t+a-b) - f(t+a) \right| + \sigma \lambda \nonumber\\
& \leq & M + |b| + \lambda.
\end{eqnarray}
Then, by (\ref{June2_1}) and (\ref{Est1}) of Lemma \ref{LEMEst1} with $\delta=-\sigma \lambda$, we have
\[
\frac{a^2+b^2}{2} -  \frac{\sigma \lambda (a+b)}{2} + \frac{\lambda^2}{4} \leq E \leq \frac{a^2-b^2}{2} +  \nu^2,
\]
which gives by (\ref{20May5}) and (\ref{20May3})
\begin{eqnarray}\label{19MayEq1}
b^2 & \leq  & \nu^2 - \frac{\lambda^2}{4} + \frac{\sigma \lambda (a+b)}{2}\nonumber\\
& \leq & \nu^2  - \frac{\lambda^2}{4} + \frac{\lambda (a+|b|)}{2} \nonumber\\
& \leq & \nu^2  - \frac{\lambda^2}{4}  + \frac{\lambda M}{2} + \frac{\lambda^2}{2} + \lambda |b| \nonumber\\
& \leq & \nu^2  + \frac{\mu^2}{4} + \frac{\mu M}{2} + \mu |b|.
\end{eqnarray}
The roots of the quadratic polynomial $b^2 -\mu b - \nu^2 - \mu^2/4 - \mu M/2$ (in the $b$ variable) are given by 
\[
\frac{\mu - \sqrt{2 \mu^2+2\mu M + 4 \nu^2}}{2} \quad \mbox{and} \quad \frac{\mu + \sqrt{2 \mu^2+2\mu M + 4 \nu^2}}{2}
\]
so the inequality (\ref{19MayEq1}) implies that
\begin{multline}
|b| \leq \frac{\mu + \sqrt{2 \mu^2+2\mu M + 4 \nu^2}}{2} = \rho M D \\
\, \mbox{ with } \quad D=D(\rho,B) := B\rho + \sqrt{2B^2\rho^2 + B+1}.
\end{multline}
By $1$-lipschitzness of $f$ and $g$ together with (\ref{20May1}), we infer that 
\begin{eqnarray}\label{June2_3}
a & = & g(t+a-b) - f(t+a) + \sigma \lambda \nonumber\\
& = & g(\tau) - f(\tau)  + f(t+a-b) - f(t+a)  + \sigma \lambda \nonumber\\
& \geq & M - 2 (\tau -T) - |b| - \mu,
\end{eqnarray}
which yields (by (\ref{20May0}))
\begin{eqnarray*}
\tau - a & \leq &\tau - M + 2 (\tau -T) + |b| + \mu \\
& \leq & T+  \frac{M}{4} +\frac{M}{4B}   - M + 2 \left( \frac{M}{4} +\frac{M}{4B} \right) + \rho MD + 2B \rho^2M \\
& = & T - \frac{M}{4} \left( 1 - \frac{3}{B} - 4 \rho D - 8 B \rho^2 \right).
\end{eqnarray*}
In conclusion, we have proved that if 
\begin{eqnarray}\label{June2_2}
1 - \frac{3}{B} - 4 \rho D - 8 B \rho^2 >0,
\end{eqnarray}
then we have $\tau -a <T < \tau$ and Lemma \ref{LEMEst2} can be applied to the $2$-Lipschitz function $h:=g-f$. Assuming that (\ref{June2_2}) holds, we obtain  
\begin{multline*}
\int_{\tau-a}^{\tau}h(t)\, dt \geq \\
\frac{h(T)^2}{4} - \frac{a^2}{2} +  \frac{h(\tau-a)^2+h(\tau)^2}{8}  - \frac{h(T) (h(\tau-a)+h(\tau))}{4} \\
+\frac{a\left( h(T) -2T+2\tau \right)}{2}  +  \frac{(T-\tau+a) h(\tau -a) }{2} + \frac{(\tau - T)h(\tau)}{2} -(\tau-T)^2,
\end{multline*}
where (remember (\ref{20May1}))
\[
\begin{array}{rcl}
h(T) & = & M\\
h(\tau) & = &  a + \delta_1 \quad \mbox{with} \quad \delta_1:=  -\sigma \lambda + f(t+a) - f(t+a-b)\\
h(\tau -a) & = & a + \delta_2 \quad \mbox{with} \quad \delta_2 :=   f(t) - f(t-b).
\end{array}
\]
So, by setting $u:=\tau -T$, we have
\begin{multline}\label{June2_4}
 \int_{ \tau -a}^{\tau} h(s) \, ds \geq 
\frac{M^2}{4} - \frac{a^2}{2} +  \frac{(a+\delta_2)^2+(a+ \delta_1)^2}{8}  - \frac{M (2a+\delta_1 + \delta_2)}{4} \\
+\frac{a\left( M +2u \right)}{2}  +  \frac{(a-u) (a+\delta_2) }{2} + \frac{u(a+\delta_1)}{2} -u^2 \\
=   \frac{a^2}{2} +  \frac{M^2-a^2}{4} + au  -u^2 \\
+  \frac{(a-M)(\delta_1+3 \delta_2)}{4} + \frac{\delta_1^2 + \delta_2^2}{8} + \frac{u (\delta_1 -  \delta_2)}{2}+ \frac{M\delta_2}{2}.
\end{multline}
By construction, we have $u\in [M/4(1-1/B),M/4(1+1/B)]$ (see (\ref{20May0})) and, by (\ref{20May3})  and (\ref{June2_3}),
\[
M (1-F) - 2 u  \leq a \leq M (1+  F)
\]
with 
\[
F=F(\rho,B) :=  \rho D + 2B\rho^2 >0.
\]
So Lemma \ref{LEMEst3} gives
\begin{eqnarray}\label{June3_3}
 \frac{M^2-a^2}{4} + a u  -u^2 \geq  \frac{M^2}{16} \left( 3 -\frac{2 }{B}  -\frac{1}{B^2} \right) - \frac{M^2F}{4} \left(1 + F + \frac{1}{B} \right).
\end{eqnarray}
We need now to bound from below the remaining terms of the right-hand side of (\ref{June2_4}). We have by $1$-lipschitzness of $f$ and $g$ and the above inequalities 
\[
|\delta_1 | \leq \lambda + \left|  f(t+a) - f(t+a-b)\right| \leq \mu + |b| \leq  M \left( 2B\rho^2 + \rho  D \right),
\]
\[
|\delta_2 | \leq  \left|  f(t) - f(t-b)\right| \leq |b|  \leq M  \rho  D,
\]
\[ 
0 \leq u \leq  \frac{M}{4} \left( 1+ \frac{1}{B} \right),
\]
and
\begin{eqnarray*}
|a-M| & = & \left| g(t+a-b) - f(t+a) + \sigma \lambda -g(T)+f(T) \right| \\
& \leq &  \left| g(t+a-b) - g(T) \right| + \left| f(T)  - f(t+a) \right| + \lambda \\
& \leq & |\tau -T| + | T - \tau -b| + \mu \\
& \leq & 2 |\tau -T| + |b| + \mu \\
& \leq & \frac{M}{2} \left( 1+ \frac{1}{B} \right)  + \rho M D +  2B\rho^2M.
\end{eqnarray*}
So we infer that
\begin{eqnarray}\label{June3_4}
& \quad &   \frac{(a-M)(\delta_1+3 \delta_2)}{4} + \frac{\delta_1^2 + \delta_2^2}{8} + \frac{u (\delta_1 -  \delta_2)}{2}+ \frac{M\delta_2}{2}\nonumber\\
& \geq & - \frac{|a-M| \left(\delta_1|+3 |\delta_2|\right)}{4} - \frac{u\left(|\delta_1|+|\delta_2|\right)}{2} -  \frac{M|\delta_2|}{2}\nonumber\\
& \geq & - \frac{M^2}{4} \left(  \frac{1}{2} \left( 1+ \frac{1}{B} \right)  + \rho D +  2B\rho^2\right)   \left( 2B\rho^2 + 4 \rho  D \right)\nonumber\\
& \quad & \qquad \qquad \qquad - \frac{M^2}{8}  \left( 1+ \frac{1}{B} \right)  \left( 2B\rho^2 + 2 \rho  D \right) - \frac{M^2}{2} \rho D.
\end{eqnarray}
Finally, we note that $E$ can be written as 
 \begin{eqnarray*}
 E & = & \int_{ t-b}^{t+a-b} g(s) -f(s)\, ds  + \int_{t-b}^{t} f (s) \, ds  - \int_{t+a-b}^{t+a } f (s) \, ds \\
 & =  & \int_{ \tau -a}^{\tau} h(s) \, ds + \int_{t-b}^{t} f (s) \, ds  - \int_{t+a-b}^{t+a } f (s) \, ds,
 \end{eqnarray*}
where by $1$-lipschitness of $f$ we have
\[
\int_{t-b}^{t} f (s) \, ds  - \int_{t+a-b}^{t+a } f (s) \, ds =  \int_{t-b}^{t} f (s) -f(s+a) \, ds \geq - |b| a,
\]
and we infer  that $E$ satisfies 
\begin{eqnarray}\label{June3_5}
E & \geq &\int_{\tau-a}^{\tau}h(t)\, dt  - |b| a,
\end{eqnarray}
where the term $\int_{\tau-a}^{\tau}h(t)\, dt$ can be bounded from below thanks to (\ref{June2_4}), (\ref{June3_3}) and (\ref{June3_4}).

In conclusion, we have proved that if  (\ref{June2_2}) is satisfied, then, by (\ref{June2_1}), (\ref{June3_5}) and the related inequalities, we have
\begin{eqnarray*}
\frac{a^2}{2}  \geq \frac{a^2-b^2}{2}  & \geq & E -\rho^2M^2\\
& \geq &\int_{\tau-a}^{\tau}h(t)\, dt  - |b| a -\rho^2M^2 \\
& \geq &\int_{\tau-a}^{\tau}h(t)\, dt  - |b| |a-M| - |b| M -\rho^2M^2\\
& \geq & \frac{a^2}{2} +  \frac{M^2}{16} \left( 3 -\frac{2 }{B}  -\frac{1}{B^2} \right) - \frac{M^2}{4} G,
\end{eqnarray*}
where (remembering that $F(\rho,B) =  \rho D + 2B\rho^2$)
\begin{multline*}
G=G(\rho,B) := F \left(1 + F + \frac{1}{B} \right) \\
+  \left(   1+ \frac{1}{B}   + 2 \rho D +  4B\rho^2\right)   \left( B\rho^2 + 2 \rho  D \right)\\
+    \left( 1+ \frac{1}{B} \right)  \left( B\rho^2 +  \rho  D \right) + 2 \rho D\\
+  2 \rho  D  \left( \left( 1+ \frac{1}{B} \right)  + 2 \rho  D +  4B\rho^2 \right) + 4 \rho  D +4 \rho^2\\
= 6 \left(2+\frac{1}{B} \right) \rho D +  4(2+B) \rho^2 + 9 \rho^2 D^2 + 22 B \rho^3 D + 8B^2 \rho^4.
\end{multline*}
and $D=B\rho + \sqrt{2B^2\rho^2 + B+1}$.
We obtain a contradiction if the pair $\rho, B \in (0,1/8)\times [2,\infty)$ satisfies 
\begin{eqnarray}\label{June3_6}
 1 - \frac{3}{B} - 4 \rho D - 8 B \rho^2 >0 \quad \mbox{and} \quad \frac{M^2}{16} \left( 3 -\frac{2 }{B}  -\frac{1}{B^2} \right) - \frac{M^2}{4} G(\rho,B) >0.
\end{eqnarray}
Since  we have for every $B> 3$,
\[
 1 - \frac{3}{B} >0, \quad 3 -\frac{2 }{B}  -\frac{1}{B^2}  >0 \quad \mbox{and} \quad \lim_{\rho \downarrow 0} 4 \rho D + 8 B \rho^2 = \lim_{\rho \downarrow 0} G(\rho,B) = 0,
\]
such pairs exist for any choice of $B$ in $(3,\infty)$. For example, if we take $B=4$, then (\ref{June3_6}) is equivalent to requiring that $\rho\in (0,1/8)$  satisfies
\[
\frac{1}{4} > 4 \rho D + 32 \rho^2
\]
and
 \begin{multline*}
 \frac{39}{16} > 4G(\rho,4) = 54 \rho D + 96 \rho^2 + 36 \rho^2 D^2 + 352 \rho^3 D + 512 \rho^4 \\
  \mbox{with}  \quad D= 4\rho + \sqrt{32 \rho^2+5}.
 \end{multline*}
 Those properties are satisfied for $\rho=0.018$. 
 
\section{Estimates}\label{SECEst}

We gather here the technical lemmas used in the proof of Theorem \ref{MainTHM}. 

\begin{lemma}\label{LEMEst1}
Let $f,g : \R \rightarrow \R$ be $1$-Lipschitz functions and $t, a, b, \delta \in \R$  such that
\begin{eqnarray}\label{HYP1}
a>0, \quad |b| \leq a, 
\end{eqnarray}
and
\begin{eqnarray}\label{HYP2}
f(t+a)=f(t)+b, \, g(t-b)= f(t)+a, \, g(t+a-b)=f(t)+a+b +\delta.
\end{eqnarray}
Then 
\begin{eqnarray}\label{Est1}
 \int_{t-b}^{t+a-b} g(s)\, ds - \int_t^{t+a} f(s)\, ds \geq \frac{a^2+b^2}{2} + \frac{\delta (a+b)}{2} + \frac{\delta^2}{4} 
 \end{eqnarray}
 and
 \begin{eqnarray}\label{Est2}
 \int_{t-b}^{t+a-b} g(s)\, ds - \int_t^{t+a} f(s)\, ds \leq \frac{3a^2-b^2}{2} + \frac{\delta (a-b)}{2} - \frac{\delta^2}{4}.
\end{eqnarray}
\end{lemma}

\begin{proof}[Proof of Lemma \ref{LEMEst1}]
We first note that if $h:\R \rightarrow \R$ is a $1$-Lipschitz function then we have for every $c,d\in \R$ with $c\leq d$, 
\begin{eqnarray}\label{PROPh}
\int_c^d h(s)\, ds \geq \frac{1}{4} \left(h(d)-h(c)\right)^2 + \frac{1}{2} (d-c) \left( h(d)+h(c)\right) -\frac{1}{4} (d-c)^2. 
\end{eqnarray}
As a matter of fact, given $c,d \in \R$ with $c\leq d$, we can define the functions $\phi_1, \phi_2: [c,d] \rightarrow \R$ by 
\[
\phi_1(s) = h(c) - (s-c) \quad \mbox{and} \quad \phi_2(s) = h(d) + (s-d),
\] 
for all $s\in \R$ and notice that since $h$ is $1$-Lipschitz and $h(c)=\phi_1(c)$, $h(d)=\phi_2(d)$, we have 
\[
h(s) \geq h(c) - \left| s-c\right| =\phi_1(s) \quad \mbox{and} \quad h(s) \geq g(d) - \left| s-d\right| = \phi_2(s) \qquad \forall s \in [c,d].
\]
Since $\phi_1$ and $\phi_2$ are affine with different slopes, there is a unique $\bar{s}\in \R$ such that $\phi_1(\bar{s})=\phi_2(\bar{s})$, it is given by
\[
\bar{s}  =  \frac{1}{2} \Bigl( h(c)-h(d)+c+d  \Bigr).
\]
Since $\bar{s}-c= (h(c)-h(d)-c+d)/2$ and $d-\bar{s}=(h(d)-h(c)-c+d)/2\geq 0$ by $1$-Lipschitzness of $h$, we have $\phi_1\geq \phi_2$ on $[c,\bar{s}]$ and $\phi_2\geq \phi_1$ on $[\bar{s},d]$. Then, we have 
\begin{eqnarray*}
\int_{c}^{d} h(s) \, ds \geq \int_{c}^{\bar{s}} \phi_1(s)\, ds + \int_{\bar{s}}^{d} \phi_2(s)\, ds,
\end{eqnarray*}
where
\begin{eqnarray*}
 \int_{c}^{\bar{s}} \phi_1(s)\, ds & = & \left( \bar{s}-c \right) \left( h(c)+c\right) - \frac{1}{2} \Bigl(\bar{s}^2 -c^2 \Bigr) \\
 & = & \left(\bar{s}-c\right) h(c) - \frac{1}{2} \left( \bar{s}-c \right)^2\\
 & = & \frac{1}{2} (h(c)-h(d)-c+d) h(c)  - \frac{1}{8} (h(c)-h(d)-c+d)^2
\end{eqnarray*}
and
\begin{eqnarray*} 
 \int_{\bar{s}}^{d} \phi_2(s)\, ds & = & (d-\bar{s}) (h(d)-d) +\frac{1}{2} \left( d^2-\bar{s}^2\right)\\
  & = &  \left(d-\bar{s} \right) h(d)    - \frac{1}{2} \left( d-\bar{s} \right)^2\\
  & = & \frac{1}{2} (h(d)-h(c)-c+d)  h(d)    - \frac{1}{8} (h(d)-h(c)-c+d)^2,
\end{eqnarray*}
which gives (\ref{PROPh}). We can apply (\ref{PROPh}) to $h=g$ on the interval $[t-b,t+a-b]$ and $h=-f$ on $[t,t+a]$. We obtain 
\begin{multline*}
\int_{t-b}^{t+a-b} g(s)\, ds - \int_t^{t+a} f(s)\, ds \geq \\
 \frac{1}{4} \left(g(t+a-b)-g(t-b)\right)^2 + \frac{a}{2}\left( g(t+a-b)+g(t-b)\right) -\frac{a^2}{4}\\
 +  \frac{1}{4} \left(-f(t+a)+f(t)\right)^2 + \frac{a}{2}\left( -f(t+a)-f(t)\right) -\frac{a^2}{4},
\end{multline*}
which, by (\ref{HYP2}), gives
\begin{multline*}
\int_{t-b}^{t+a-b} g(s)\, ds - \int_t^{t+a} f(s)\, ds \geq \\
 \frac{1}{4} \left(b+\delta \right)^2 + \frac{a}{2} \left( 2f(t)+2a+b+\delta\right)  -\frac{a^2}{4} + \frac{b^2}{4} +\frac{a}{2} \left( -2f(t) -b\right) -\frac{a^2}{4}
 \end{multline*}
 and implies  (\ref{Est1}). We can also apply  (\ref{PROPh}) to $h=-g$ on the interval $[t-b,t+a-b]$ and $h=f$ on $[t,t+a]$ to get
 \begin{multline*}
\int_{t-b}^{t+a-b} -g(s)\, ds + \int_t^{t+a} f(s)\, ds \geq \\
 \frac{1}{4} \left(-g(t+a-b)+g(t-b)\right)^2 + \frac{a}{2}\left(- g(t+a-b)-g(t-b)\right) -\frac{a^2}{4}\\
 +  \frac{1}{4} \left(f(t+a)-f(t)\right)^2 + \frac{a}{2}\left( f(t+a)+ f(t)\right) -\frac{a^2}{4}.
\end{multline*}
By (\ref{HYP2}), we obtain
\begin{multline*}
\int_{t-b}^{t+a-b} g(s)\, ds - \int_t^{t+a} f(s)\, ds \leq \\
  - \frac{1}{4} \left(b+\delta \right)^2 + \frac{a}{2} \left( 2f(t)+2a+b+\delta\right) + \frac{a^2}{4} - \frac{b^2}{4} +\frac{a}{2} \left( -2f(t) -b\right) + \frac{a^2}{4},
\end{multline*}
which  gives (\ref{Est2}). 
\end{proof}

\begin{lemma}\label{LEMEst2}
Let $h:\R\rightarrow \R$ be a $2$-Lipschitz function and $T, \tau, a\in \R$ such that 
\begin{eqnarray}\label{HYP3}
\tau -a < T < \tau.
\end{eqnarray}
Then we have 
\begin{multline}\label{Est3}
\int_{\tau-a}^{\tau}h(t)\, dt \geq \\
\frac{h(T)^2}{4} - \frac{a^2}{2} +  \frac{h(\tau-a)^2+h(\tau)^2}{8}  - \frac{h(T) (h(\tau-a)+h(\tau))}{4} \\
+\frac{a\left( h(T) -2T+2\tau \right)}{2}  +  \frac{(T-\tau+a) h(\tau -a) }{2} + \frac{(\tau - T)h(\tau)}{2} -(\tau-T)^2.
\end{multline}
\end{lemma}

\begin{proof}[Proof of Lemma \ref{LEMEst2}]
Since the function $h/2$ is $1$-Lipschitz, we can apply the lower bound (\ref{PROPh}) obtained at the beginning of the proof of Lemma \ref{LEMEst1}. We obtain that for every $c,d\in \R$, with $c\leq d$, we have 
\begin{eqnarray}\label{PROPh2}
\int_c^d h(s)\, ds \geq \frac{1}{8} \left(h(d)-h(c)\right)^2 + \frac{1}{2} (d-c) \left( h(d)+h(c)\right) -\frac{1}{2} (d-c)^2. 
\end{eqnarray}
We infer that
\begin{eqnarray*}
& \quad & \int_{\tau-a}^{\tau} h(t)\, dt \\
& = &\int_{\tau-a}^T h(t)\, dt + \int_{T}^{\tau} h(t)\, dt \\
& \geq &  \frac{1}{8} \left(h(T)-h(\tau -a)\right)^2 + \frac{1}{2} (T-\tau +a) \left( h(T)+h(\tau -a)\right) -\frac{1}{2} (T-\tau +a)^2 \\
& \quad & + \frac{1}{8} \left(h(\tau)-h(T)\right)^2 + \frac{1}{2} (\tau-T) \left( h(\tau)+h(T)\right) -\frac{1}{2} (\tau -T)^2
\end{eqnarray*}
which gives (\ref{Est3}).
\end{proof}

\begin{lemma}\label{LEMEst3}
Let $M, \delta>0$ and $B>1$ be fixed. Then for any $a,u \in \R$ such that 
\[
M (1-\delta) -2u \leq a \leq M (1+\delta) \quad \mbox{and} \quad \frac{M}{4}\left( 1-\frac{1}{B} \right)   \leq u\leq \frac{M}{4} \left( 1+\frac{1}{B} \right) 
\]
we have
\begin{eqnarray}
 \frac{M^2-a^2}{4} + a u  -u^2 \geq  \frac{M^2}{16}\left( 3-\frac{2}{B}-\frac{1}{B^2}\right) - \frac{M^2\delta}{4} \left( 1 + \delta + \frac{1}{B}  \right).
\end{eqnarray}
\end{lemma}

\begin{proof}[Proof of Lemma \ref{LEMEst3}]
Let $\Phi : \R^2 \rightarrow \R$ be the function defined by 
\[
\Phi(a,u) :=  \frac{M^2-a^2}{4} + u (a -u)  = \frac{M^2}{4} - \frac{ ( a -2u)^2}{4} \qquad \forall (a,u) \in \R^2.
\] 
The gradient of $\Phi$ is given by
\[
\nabla \Phi (a,u) = - \frac{\left( a -2u \right)}{2} \left( \begin{matrix} 1 \\ -2 \end{matrix} \right) 
\]
So, if it vanishes at $(a,u)$, we have $\Phi(a,u)=M^2/4$ and $(a,u)$ is necessarily a local maximum of $\phi$. As a consequence, the minimum of $\Phi$ on the closed set 
\[
F = \left\{ (a,u) \in \R^2 \, \vert \, M(1-\delta) -2u \leq a \leq M (1+  \delta),  \frac{M}{4} \left( 1-\frac{1}{B} \right) \leq u\leq \frac{M}{4} \left( 1+\frac{1}{B} \right) \right\}
\] 
is attained on the boundary. We check easily that the vector $(1,-2)$ is never orthogonal to the fours faces of the boundary. This shows that the minimum of $\Phi$ on $F$ has to be attained at one of the four corners of $F$ whose images by $\Phi$ are given by 
\begin{eqnarray*}
& \quad &  \Phi\left(M(1+ \delta), \frac{M}{4} \left(1+ \frac{\epsilon}{B} \right) \right) \\
 & = & \frac{M^2}{4} \left( - 2\delta - \delta^2\right) +   \frac{M}{4} \left(1+ \frac{\epsilon}{B} \right)  \left( M (1+ \delta) - \frac{M}{4}  \left(1+ \frac{\epsilon}{B} \right) \right) \\
 & = & \frac{M^2}{16} \left( 3+\frac{2 \epsilon}{B} - \frac{1}{B^2}  \right) + \frac{M^2}{4} \left( -\delta - \delta^2 + \frac{\epsilon \delta}{B} \right)
\end{eqnarray*}
 and
\begin{eqnarray*}
& \quad &  \Phi\left(M(1- \delta) -  \frac{M}{2} \left(1+ \frac{\epsilon}{B} \right) , \frac{M}{4} \left(1+ \frac{\epsilon}{B} \right) \right) \\
 & = & \frac{M^2}{4} - \frac{1}{4} \left( M (1- \delta) - M \left(1+ \frac{\epsilon}{B}   \right)   \right)^2 \\
 & = & \frac{M^2}{4} \left( 1-\frac{1}{B^2} \right) - \frac{M^2}{4} \left(\frac{\epsilon \delta}{B} + \delta^2 \right)
\end{eqnarray*}
for $\epsilon = \pm 1$. We note that since $B>1$, we have $4(1-1/B^2)>3-2/B-1/B^2$, so the minimum of the four values above is attained for $(a,u) = (M(1+\delta),M(1-1/B)/4)$, so we obtain that 
\begin{eqnarray*}
\Phi(a,u)  \geq     \frac{M^2}{16}  \left( 3+\frac{2 \epsilon}{B} - \frac{1}{B^2}  \right)     - \frac{M^2\delta}{4} \left(  \delta + \delta^2 + \frac{\delta}{B} \right),
\end{eqnarray*}
which proves the result.
\end{proof}


\begin{thebibliography}{99}

\bibitem{aa18}
A.~Akopyan and S.~Avvakumov.
\newblock {\em Any cyclic quadrilateral can be inscribed in any closed convex smooth curve}.
\newblock Forum Math. Sigma, 6, 2018.

\bibitem{gg73}
M.~Golubitsky and V.~Guillemin.
\newblock Stable mapping and their singularities. 
\newblock Graduate Texts in Mathematics, Vol. 14. Springer-Verlag, New-York-Heidelberg, 1973. 

\bibitem{karasev13}
R.N.~Karasev.
\newblock {\em On two conjectures of Makeev}. 
\newblock Translated from Zap. Nauchn. Sem. S.-Peterburg. Otdel. Mat. Inst. Steklov. (POMI) 415 (2013), Geometriya i Topologiya. 12, 5–14 J. Math. Sci.
(N.Y.) 212(5): 521–-526, 2016.

\bibitem{matschke14}
B.~Matschke.
\newblock {\em A survey on the square peg problem}.
\newblock Notices Amer. Math. Soc., 61(4):346--352, 2014.

\bibitem{matschke20}
B.~Matschke.
\newblock {\em Quadrilaterals inscribed in convex curves}.
\newblock Preprint, https://arxiv.org/abs/1801.01945, 2020.

\bibitem{tao17}
T.~Tao.
\newblock  {\em An integration approach to the Toeplitz square peg problem}.
\newblock Forum Math. Sigma, 5, 2017.

\bibitem{toeplitz}
O.~Toeplitz.
\newblock {\em \"{U}ber einige aufgaben der Analysis situs}.
\newblock Verhandlungen der Schweizerischen Naturforschenden Gesellschaft in Solothurn, 4:197, 1911.
\end{thebibliography}
\end{document}